\documentclass[12pt,a4paper,twoside,leqno]{article}
\usepackage{epsfig, cite}
\usepackage{graphpap,latexsym,epsf}
\usepackage{graphics}
\usepackage[usenames,dvipsnames]{color}
\def\betti #1{{\color{red}#1}}
\def\juerg #1{{\color{green}#1}}
\def\pier #1{{\color{blue}#1}}
\def\bettinew #1{{\color{magenta}#1}}
\let\betti\relax
\let\pier\relax
\let\juerg\relax
\let\bettinew\relax
\usepackage{amsthm,amsmath,amssymb}

\definecolor{green1}{rgb}{0.1,0.65,0}
\definecolor{blue1}{rgb}{0.14,0.6,1.0}
\definecolor{blue3}{rgb}{0.05,0.05,0.5}
\definecolor{whit}{rgb}{1,1,1}
\definecolor{grey1}{rgb}{0.7,0.7,0.7}

\def\trait #1 #2 #3 {\vrule width #1pt height #2pt depth #3pt}
\def\fin{
      \trait .3 5 0
      \trait 5 .3 0
      \kern-5pt
      \trait 5 5 -4.7
      \trait 0.3 5 0
\medskip}
\newtheorem{teor}{Theorem}[section]
\newtheorem{defin}[teor]{Definition}
\newtheorem{lemm}[teor]{Lemma}
\newtheorem{osse}[teor]{Remark}
\newtheorem{prop}[teor]{Proposition}
\newtheorem{defi}[teor]{Definition}
\newtheorem{coro}[teor]{Corollary}
\newtheorem{prob}[teor]{Problem}
\newtheorem{hypo}[teor]{Hypothesis}

\newcommand{\bele}{\begin{lemm}\begin{sl}}
\newcommand{\enle}{\end{sl}\end{lemm}}
\newcommand{\bedef}{\begin{defi}\begin{sl}}
\newcommand{\eddef}{\end{sl}\end{defi}}
\newcommand{\bete}{\begin{teor}\begin{sl}}
\newcommand{\ente}{\end{sl}\end{teor}}
\newcommand{\beos}{\begin{osse}\begin{rm}}
\newcommand{\eddos}{\end{rm}\end{osse}}
\newcommand{\bepr}{\begin{prop}\begin{sl}}
\newcommand{\empr}{\end{sl}\end{prop}}
\newcommand{\bepro}{\begin{prob}\begin{rm}}
\newcommand{\empro}{\end{rm}\end{prob}}
\newcommand{\bede}{\begin{defin}\begin{sl}}
\newcommand{\edde}{\end{sl}\end{defin}}
\newcommand{\beco}{\begin{coro}\begin{sl}}
\newcommand{\enco}{\end{sl}\end{coro}}
\newcommand{\behy}{\begin{hypo}\begin{sl}}
\newcommand{\enhy}{\end{sl}\end{hypo}}
\newcommand{\beq}{\begin{equation}} 
\newcommand{\eeq}{\end{equation}}
\newcommand{\beqa}{\begin{eqnarray}}
\newcommand{\eeqa}{\end{eqnarray}}

\newcommand{\thspace}{\hspace{3mm}}

\newcommand{\RR}{\mathbb{R}}

\newcommand{\NN}{\mathbb{N}}

\newcommand{\dx}{\,{\rm d}x}
\newcommand{\dt}{\,{\rm d}t}
\newcommand{\ds}{\,{\rm d}s}

\newcommand{\dz}{\,{\rm d}z}

\def\qed{\ifmmode 
  \else \leavevmode\unskip\penalty9999 \hbox{}\nobreak\hfill
  \fi
  \quad\hbox{\hskip.5em\vrule width.4em height.6em depth.05em\hskip.1em}}
\def\endproofsym{\qed}

\def\endnobox{\def\endproofsym{}\end{proof}\def\endproofsym{\qed}}

\newcommand{\no}{\nonumber}

\newcommand{\dn}{\partial_{\bf n}}
\newcommand{\doma}{\partial\Omega}
\newcommand{\oma}{\Omega}

\newcommand{\ioma}{\int_\Omega}

\newcommand{\iTT}{\int_0^T}
\newcommand{\txinto}{\int_0^t\!\!\ioma}

\newcommand{\sumk}{\sum_{k=1}^n}

\newcommand{\lzo}{L^2(\Omega)}
\newcommand{\lzq}{L^2(Q)}
\newcommand{\lio}{L^\infty(\Omega)}
\newcommand{\liq}{L^\infty(Q)}
\newcommand{\heins}{H^1(\Omega)}
\newcommand{\hzwei}{H^2(\Omega)}
\newcommand{\lzv}{L^2(0,T;V)}
\newcommand{\lzvt}{L^2(0,t;V)}
\newcommand{\lzh}{L^2(0,T;H)}
\newcommand{\lzht}{L^2(0,t;H)}

\newcommand{\liv}{L^\infty(0,T;V)}

\newcommand{\lzhz}{L^2(0,T;H^2(\Omega))}

\newcommand{\pt}{\partial_t}

\newcommand{\s}{\sigma}

\newcommand{\PP}{\bf {(CP)}}
\newcommand{\ur}{{\cal U}_R}
\newcommand{\uad}{{\cal U}_{\rm ad}}
\newcommand{\bu}{\overline{u}}
\newcommand{\bmu}{\overline{\mu}}
\newcommand{\bs}{\overline{\sigma}}
\newcommand{\bphi}{\overline{\varphi}}

\newcommand{\bq}{\overline{Q}}
\newcommand{\nf}{\bf {n}}
\newcommand{\cs}{{\cal S}}

\let\TeXchi\chi
\def\chi{{\setbox0 \hbox{\mathsurround0pt
$\TeXchi$}\hbox{\raise\dp0 \copy0 }}}

\newcommand{\as}{({\bf H1})\mbox{--}({\bf H6})}

\def\vp{\varphi}

\newfont{\ctv}{msam10}

\def\fine{\hfill\kern4pt \vrule height4pt depth0pt width4pt }

\numberwithin{equation}{section}

\begin{document}

\title{\bf {Optimal distributed control of a diffuse interface model of tumor 
growth\footnote{\pier{The financial support of the FP7-IDEAS-ERC-StG \#256872
(EntroPhase) is gratefully acknowledged. The paper 
also benefits from the support of the MIUR-PRIN Grant 2010A2TFX2 ``Calculus of Variations'' for PC and GG, and the GNAMPA (Gruppo Nazionale per l'Analisi Matematica, la Probabilit\`a e le loro Applicazioni) of INdAM (Istituto Nazionale di Alta Matematica) for PC, GG and ER.}}}}

\date{}

\author{
Pierluigi Colli\footnote{Dipartimento di Matematica ``F. Casorati,'' Universit\`{a}
di Pavia, and IMATI-C.N.R., Via Ferrata 1, I-27100 Pavia, Italy (\tt pierluigi.colli@unipv.it).} , 
Gianni Gilardi\footnote{Dipartimento di Matematica ``F. Casorati,'' Universit\`{a}
di Pavia, and IMATI-C.N.R., Via Ferrata 1, I-27100 Pavia, Italy (\tt gianni.gilardi@unipv.it).} ,
Elisabetta Rocca\footnote{Weierstrass Institute for Applied
Analysis and Stochastics, Mohrenstrasse~39, D-10117 Berlin,
Germany ({\tt  rocca@wias-berlin.de}), and Dipartimento di Matematica, 
Universit\`{a} di Milano,
Via Saldini 50, I-20133 Milano, Italy ({\tt elisabetta.rocca@unimi.it}).}\\
and J\"urgen Sprekels\footnote{Weierstrass Institute for Applied
Analysis and Stochastics, Mohrenstrasse~39, D-10117 Berlin,
Germany ({\tt  sprekels@wias-berlin.de}), and Institut f\"ur Mathematik der Humboldt-Universit\"at
zu Berlin, Unter den Linden 6, D-10099 Berlin, Germany\pier{.}}}

\maketitle

\vspace{-9mm}

\noindent {\bf Abstract.} 
In this paper, a  distributed optimal control problem is studied for a diffuse interface model of tumor growth
which was proposed by Hawkins--Daruud et al. in \cite{HZO}. The model consists of a Cahn--Hilliard equation for the tumor
cell fraction $\vp$ coupled to a reaction-diffusion equation for a function $\s$ representing the nutrient-rich 
extracellular water volume fraction. The distributed control $u$ monitors as a right-hand side the equation for
$\s$ and can be interpreted as a nutrient supply or a medication, while the cost function, which is of standard 
tracking type, is meant to keep the tumor cell fraction under control during the evolution. 
We show that the control-to-state operator is Fr\'echet differentiable between appropriate Banach spaces and derive the first-order necessary optimality conditions in terms of a variational inequality involving the adjoint state 
variables.

\vspace{.4cm}

\vspace{4mm}

\noindent
{\bf AMS (MOS) subject clas\-si\-fi\-ca\-tion:}\thspace
\pier{35K61,49J20, 49K20, 92C50.}

\pagestyle{myheadings}
\newcommand\testopari{\sc \pier{Colli \ --- \ Gilardi \ --- \ Rocca \ --- \ Sprekels}}
\newcommand\testodispari{\sc \pier{Control of a diffuse interface model of tumor growth}}
\markboth{\testodispari}{\testopari}

\section{Introduction}\label{intro}
Let $\oma\subset\RR^3$ be an open bounded and connected set with a smooth boundary $\doma$, and let 
$\nf$ denote the outward unit normal to $\doma$. Moreover, let a fixed final time $T>0$ be given and 
$Q:=\oma\times (0,T)$, $\Sigma:=\doma\times (0,T)$. We investigate in this paper the following distributed optimal control problem:

\vspace{2mm}
\noindent$\PP$ \,\,\,Minimize the cost functional
\begin{align}
\label{cost}
{\cal J}(\vp,u)=&\ \frac{\beta_Q}2\iTT\!\!\ioma |\vp-\vp_Q|^2\dx\dt\,+\,
\frac{\beta_\Omega}2\ioma |\vp(T)-\vp_\oma|^2\dx\\
&+\,\frac{\beta_u}2\iTT\ioma |u|^2\dx\dt\nonumber 
\end{align}
subject to the control constraint
\begin{equation} \label{Uad}
u\in\uad:=\{u\in L^\infty(Q):\,\,u_{\rm {min}}\le u\le u_{\rm max} \,\,\,\mbox{a.\,e. in }\,Q\}
\end{equation}
and to the state system
\begin{align}
\label{ss1}
&\vp_t-\Delta\mu= P(\vp)(\s-\delta\mu) \quad\mbox{in }\,Q,\\[1mm]
\label{ss2}
&\mu=-\Delta\vp+F'(\vp)\quad\mbox{in }\,Q,\\[1mm]
\label{ss3}
&\s_t-\Delta \s=-P(\vp)(\s-\delta\mu)\betti{{}+u{}} \quad\mbox{in }\,Q,\\[1mm]
\label{ss4}
&\partial_{\nf}\vp=\partial_{\nf}\mu=\partial_{\nf}\s=0 \quad\mbox{on }\,\Sigma,\\[1mm]
\label{ss5}
&\vp(0)=\vp_0,\quad\s(0) =\s_0 \quad\mbox{in }\,\oma.
\end{align} 

\noindent The quantities occurring in the above expressions have the following meaning: $\beta_Q$, $\beta_\oma$,
$\beta_u$ are nonnegative constants, $\delta \bettinew{{}>{}} 0$ is a constant, $\vp_Q\in L^2(Q)$, $\vp_\oma\in L^2(\oma)$, $u_{\rm min}
\in L^\infty(\pier{Q})$, $u_{\rm max} \in L^\infty(\pier{Q})$ are given functions \pier{such that 
$u_{\rm min} \leq u_{\rm max} $ almost everywhere in $Q$}, $F$ and $P$ are given nonlinearities,
$\dn$ denotes the derivative in the direction of the outward unit normal $\nf$, and $\vp_0$, $\s_0$
are given initial data. In the following, we will always assume that $\delta=1$, which has no bearing on the
mathematical analysis.

The state equations \eqref{ss1}--\eqref{ss5} constitute an approximation to a model for the 
dynamics of tumor growth due to \cite{HZO} (see also \cite{Hil,WZZ}) in which the velocities are set to
zero and the only state variables considered are the tumor cell fraction $\vp$ and the nutrient-rich
extracellular water fraction $\s$. Typically, the function $F$ occurring in the chemical potential $\mu$ is a double-well potential, and $P$ denotes a suitable proliferation function\betti{, which is in general a} \juerg{nonnegative} \betti{and regular function of $\vp$}.

Altogether, the optimal control problem $\PP$ can be interpreted as the search for a strategy how to apply
a control $u$ (which may represent the supply of a nutrient (see \cite{BC}), or even a drug in a  chemotherapy)
properly in order that
\begin{enumerate}
\item[(i)] \pier{the integral over the full space-time domain of the squared} amount of nutrient or drug supplied (which is restricted by the control constraints) does not
inflict any harm on the patient (which is expressed by the presence of the third summand in the cost
functional)\pier{;}
\item[(ii)] \betti{a desired} evolution and final distribution of the tumor cells (which is expressed by the target
functions $\vp_Q$ and $\vp_\oma$) is realized in the best possible way.
\end{enumerate}
The ratios $\beta_Q/\beta_u$ and $\beta_\oma/\beta_u$ indicate which importance the conflicting targets
`avoid unnecessary harm to the patient' and `quality of the approximation of $\vp_Q$, $\vp_\oma$' are given in the
strategy. We remark that in practice it would be safer for the patient (and thus more desirable) to 
approximate the target functions rather in the $L^\infty$ sense than in the $L^2$ sense; however, 
in view of the analytical difficulties that are inherent in the state system, this presently seems to be out 
of reach.  \pier{Of course, other integral terms depending on $\sigma$ and analogous to the ones acting on $\varphi$ could be added to \juerg{the control functional, and it is our opinion that such} a modified functional may be tractable from a mathematical point of view. However, since the problem is already quite involved and it is not clear whether the modification is really worth to be considered for applications, we prefer not to include the extra terms in the cost functional.} 

The mathematical modeling of tumor growth dynamics has drawn much attention in the past decade  (cf., e.\,g., 
\cite{AMcE,CriLow,LowFetal,OPH}). In particular, models based on continuum mixture theory have been
derived (see \cite{WLow, CBC, OHP, CLLowW, FJC, HP}), which usually lead to Cahn--Hilliard systems involving transport and reaction terms that
govern various types of concentrations, where, in particular, the reaction terms depend on the nutrient
concentration.

While there exist quite a number of numerical simulations of diffuse-interface models of tumor growth (cf., e.\,g.,
\cite[Chap.\,8]{CriLow},\betti{\cite{CLLowW, HZO, WLow, WZZ}}, and the references given 
there), there are still only a few contributions to the
mathematical analysis of the models. The first contributions in this direction dealt with the case where the
nutrient is neglected, which then leads to the so-called Cahn--Hilliard--Hele--Shaw system 
(see\betti{\cite{LowTZ,JWZ,
BCG, WW, WZ}}). \betti{Only very recently\pier{,} in the paper \cite{DFRSS}\pier{,} the authors proved existence of weak solutions 
and some rigorous sharp interface limit for a model introduced in \cite{cwsl} (cf. also \cite{CriLow, 
CLLowW, FJC, LowFetal, WLow}), where both velocities (satisfying 
a Darcy law with Korteweg term) and 
multispecies tumor fractions, as well as the nutrient evolutions are taken into account. 
\pier{Let us also quote}  the paper \cite{GetAl}, where a new model for tumor growth including different densities is introduced and a formal sharp interface limit is performed\pier{; moreover, the well-posedness of the 
related diffuse interface model is discussed in \cite{GaLa}.
Finally, in the contribution}} \cite{FGR}, the system \eqref{ss1}--\eqref{ss5}, which constitutes
an approximation of the model introduced in \cite{HZO}, was rigorously analyzed concerning well-posedness,
regularity, and asymptotic behavior. We also refer to the recent papers
\cite{CGH,CGRS1,CGRS2}, in which various `viscous' approximations of the state system have been
studied analytically.    

In this paper, we focus on the control aspect. While there exist many contributions concerning the
well-posedness of various types of Cahn--Hilliard systems, only a few deal with their
optimal control. In this connection, we mention the papers  
\cite{wn99,hw,CGS1,ColliGilSpr}, which deal with zero Neumann boundary conditions like \eqref{ss4}, 
while in the recent papers \cite{CGS1, CGS2, CFGS1, CFGS2} dynamic boundary conditions have been
studied. A number of papers also investigates optimal control problems for 
convective Cahn--Hilliard systems (cf.\betti{\cite{ZL1, ZL2, RS}}) and Cahn--Hilliard--Navier--Stokes systems
(cf. \cite{HW1, HW2, HW3, FRS}). 
\bettinew{Regarding the problem of optimal control in tumor growth models, we can quote the papers \cite{LS}, where the problem of minimizing the volume of tumor under isoperimetric contraints is considered, and \cite{BAD}, where an advection-reaction-diffusion system for leukemia development is studied.}
However, to the authors' best knowledge, optimal control problems for the system \eqref{ss1}--\eqref{ss5}
have never been studied before. 

\betti{Indeed, the main mathematical difficulties are  related to the proofs of suitable 
stability estimates of higher order (with respect to the ones already present in \cite{FGR}), which are necessary in order to prove the differentiability (in suitable spaces) 
of the control-to-state mapping. The presence of the two nonlinearities $F$ and $P$ is indeed the main challenge in the analysis. 
\juerg{Moreover, due} to the dependence on the $L^2(\Omega)$\,--\,target $\vp_\Omega$  in the final condition for the variable $p$, \juerg{which is} related to the tumor phase $\vp$, we only get 
existence for the adjoint system in \juerg{the sense that} the $p$-equation has to be intended in a weak form, mainly in the dual of the Sobolev space 
$H^2(\Omega)$ (cf. Section~\ref{control} for \juerg{further} comments on 
\juerg{this} subject). Let us finally \pier{point out} that, 
\pier{with a view to applications}, it would be 
worth \juerg{analyzing} the case of an $L^\infty$\,--\,type control functional rather than the $L^2$\,--\,one \juerg{tackled here; however,} this would bring further difficulties in solving the adjoint system \juerg{in which measures would occur on} the \pier{right-hand sides of} the equations. Hence, since the present contribution is the first one on \pier{the} control theory for 
diffuse models of tumor growth, we prefer to start with the $L^2$\,--\,control function 
$\mathcal{J}$ in \eqref{cost}. } 

\paragraph{Plan of the paper.} The paper is organized as follows: in Section 2, we formulate
the general \betti{hypotheses} and 
improve known results regarding the well-posedness and regularity of the state \betti{system}
(\ref{ss1})--(\ref{ss5}) (Theorem 2.1). We also prove a continuous dependence result 
(Theorem 2.2) which is needed for the analysis of the control problem. In  Section 3,
we study the differentiability properties of the control-to-state operator. The main 
results of this paper concerning existence and first-order necessary 
optimality conditions for the optimal control problem {\bf (CP)} are shown in Section 4. 

Throughout this paper, \pier{for a (real) Banach space $X$ we denote by $\,\|\cdot\|_X$ its norm, by $X'$ its dual space, 
and by $\langle\cdot,\cdot\rangle_X$ the dual pairing between $X'$ and $X$. 
If $X$ is an inner product space, 
then the inner product is denoted by $(\cdot,\cdot)_X$. The only exception from this convention is given
by the $L^p$ spaces, $1\le p\le\infty$, for which we use the abbreviating notation
$\|\cdot\|_p$ for the norms in both $L^p(\Omega)$ and $L^p(Q)$. Moreover, we will use the notations}
$$H:=\lzo, \quad V:=\heins, \quad W:=\{w\in \hzwei:\,\dn w=0\mbox{ a.\,e. on }\doma\}.$$
We have the dense and continuous embeddings 
$W\subset V\subset H\cong H'\subset V'\subset W'$, where $\langle u,v\rangle_V=(u,v)_H$ and
$\langle u,w\rangle_W=(u,w)_H$ for all $u\in H$, $v\in V$, and $w\in W$.

During the course of this paper, 
we will make repeated use of Young's inequality
\begin{equation}
\label{young}
a\,b\le \delta \,a^2 + \frac 1{4\delta} b^2 \quad\mbox{for all }\,
a,b\in \RR \,\mbox{ and }\, \delta>0,
\end{equation}
as well as of the fact that for three dimensions of space and smooth domains 
the embeddings $\,V\subset L^p(\Omega)$, $1\le p\le 6$, and 
$\,H^2(\Omega)\subset C^0(\overline{\Omega})$ are continuous and 
(in the first case only for $1\le p<6$) 
compact. In particular, there are positive constants $\widetilde K_i$, $i=1,2,3$,
 which depend only on the domain $\oma$, such that
\beqa\label{embed1}
\|v\|_{6}&\!\!\le\!\!&\widetilde K_1\,\|v\|_{V}\quad\forall\,
v\in V,\\[1mm]
\label{embed2}
\|v\,w\|_H&\!\!\le\!\!&\|v\|_{6}\,\|w\|_{3}\,\le\,
\widetilde K_2\,\|v\|_V\,\|w\|_V\quad\forall\,v,w\in V, 
\\[1mm]
\label{embed3}
\|v\|_{L^\infty(\oma)}&\!\!\le\!\!&\widetilde K_3\,\|v\|_{H^2(\oma)}\quad\forall\,
v\in \hzwei.
\eeqa
Moreover, we have
\beq\label{gianni}
\|v\,w\|_{V'}\,\le\,\|v\|_{W^{1,\infty}(\oma)}\,\|w\|_{V'}, \quad\forall \,v
\in W^{1,\infty}(\oma), \,\,\,\forall\,w\in V'.
\eeq
Finally, we recall that for smooth and bounded three-dimensional domains there holds the special
Gagliardo--Nirenberg inequality \betti{(cf. \cite[p. 125]{GN})}
\beq
\label{GN1}
\|v\|_3\,\le\,\widetilde{K}_4\left(\|v\|_H^{1/2}\,\|v\|_V^{1/2} \,+\,\|v\|_H\right)
\quad\,\,\forall\,v\in V\,,
\eeq 
where the positive constant $\widetilde{K}_4$ depends only on $\Omega$.

\section{General assumptions and preliminary results on the state system}
In the following, we study the state system (\ref{ss1})--(\ref{ss5}).
Since it will be convenient to rewrite various partial differential equations in this paper
as abstract equations in the framework of the Hilbert triple $(V,H,V')$,
we introduce the Riesz isomorphism $A:V\to V'$ associated with the standard scalar product of~$V$, that is,
\begin{equation}
  \langle A u, v\rangle_V= (u,v)_V
  = \int_\Omega \left(\nabla u \cdot \nabla v + uv \right)\dx
  \quad \hbox{for $u,v\in V$}.
  \label{defA}
\end{equation}
We note that the restriction of $A$ to $W$, which is given by $A u=-\Delta u+u$, for $u\in W$,
is an isomorphism from $W$ onto~$H$. Moreover, the linear operator $A$ can be continuously extended
to a linear mapping from $H$ into $W'$, where $\langle Au,v\rangle_W=(u,Av)_H$ for all $u\in H$ and
$\juerg{v}\in W$.
We also remark that, \betti{for some}  positive constant $\widetilde K_5$ which depends only on $\oma$, \betti{we have}
\begin{eqnarray}
\label{ruleA1}
\langle A u , A^{-1} v^*\rangle_V
&\!\! =\!\!&\langle v^*, u\rangle_V
  \quad \hbox{for all $u\in V$ and $v^*\in V'$,}
  \\[1mm]
\nonumber
\langle u^* , A^{-1} v^* \rangle_V
&\!\!=\!\!& (u^*,v^*)_{V'}
  \quad \hbox{for all $u^*,v^*\in V'$}, \\[1mm]
\nonumber
\|A^{-1} u^*\|_{V'} 
&\!\!\le\!\!&\widetilde K_5\,\|u^*\|_{V'} \quad\mbox{for all $u^*\in V'$}, 
\end{eqnarray}
where $(\cdot,\cdot)_{V'}$ is the dual scalar product in $V'$
associated with the standard one in~$V$. We also have, for every $v^*\in H^1(0,T; V')$,
\beq\label{ruleA2}
  \frac d{dt} \|v^*(t)\|_{V'}^2\,
  =\, 2\, \langle\partial_tv^*(t) , A^{-1} v^*(t) \rangle_V
  \quad \hbox{for almost every $t\in (0,T)$} .
\eeq

\betti{We} make for the remainder of this paper the following general assumptions on the data of the 
control problem $\PP$:

\vspace{2mm}\noindent
{\bf (H1)} \,\,\,$\beta_Q,\beta_\oma,\beta_u$ are nonnegative but not all zero.\\[2mm]
{\bf (H2)} \,\,\,$\vp_Q\in \lzq$, $\vp_\oma\in\lzo$, $u_{\rm min}\in \liq$, $u_{\rm max}\in\liq$, \,with\, 
$u_{\rm min}\le u_{\rm max}$\hfill\\
\hspace*{12.2mm} a.\,e. in $Q$.\\[2mm]
{\bf (H3)} \,\,\,$\vp_0\in H^3(\oma)$, $\s_0\in \heins$.\\[2mm]
{\bf (H4)} \,\,\,$P\in C_{\rm loc}^{1,1}(\RR)$ is nonnegative and satisfies, for almost every $s\in\RR$,
\beq
\label{2.4} 
|P'(s)|\,\le\,\alpha_1\left(1+|s|^{q-1}\right), \,\mbox{ with some $\alpha_1>0$ and some $q\in [1,4]$}.\\
\eeq
{\bf (H5)} \,\,\,$F\in C^4(\RR)$ can \,be \,written \,in \,the form \,$\,F=F_0+F_1$, where $F_0,F_1\in C^4(\RR)$,\\
\hspace*{12.2mm} and where there are constants $\alpha_i>0$, $\betti{2}\le i\le 6$, and $\rho\in [2,6)$ such that 
\beqa
\label{2.5}
&&|F_1''(s)|\le\alpha_2 \quad\forall\,s\in\RR, \\[1mm]
\label{2.6}
&&\alpha_3\left(1+|s|^{\rho-2}\right)\le F_0''(s)\le \alpha_4\left(1+|s|^{\rho-2}\right) 
\quad\forall\,s\in\RR,\\[1mm]
\label{2.7}
&&F(s)\ge \alpha_5|s| - \alpha_6 \quad\forall\,s\in\RR.
\eeqa

\noindent
The conditions ({\bf H3})--({\bf H5}) originate from the paper \cite{FGR}, where they
were postulated to guarantee the validity of some well-posedness results. We remark that not all of them are needed
for some of the results proved in \cite{FGR}; however, they are indispensable for the analysis of the control  
problem $\PP$ on which we focus in this paper. 
The following hypothesis is rather a denotation than an assumption: 

\vspace{2mm}
\noindent
{\bf (H6)} \,\,\,${\cal U}_R$ is an open set in $\lzq$ such that $\uad\subset\ur$ and $\|u\|_{\lzq}\le R$ for all
$u\in\ur$.

\vspace{3mm}
We have the following well-posedness result for the state system \eqref{ss1}--\eqref{ss5}.

\vspace{5mm} \noindent
{\sc Theorem 2.1} \quad {\em Suppose that the hypotheses} $\as$ {\em are fulfilled. Then the follwing results
hold true:}

\noindent
(i)\,\,\,\,{\em For every $u\in\ur$, the state system} \eqref{ss1}--\eqref{ss5} {\em has a unique strong solution 
triple $(\vp,\mu,\s)$ such that}
\begin{align}
\label{reguss}
&\vp\in H^1(0,T;V)\cap L^\infty(0,T;\pier{W\cap{}}H^3(\oma)), \quad \pier{\Delta \vp\in L^2(0,T;W)}, \\
&\nonumber \mu\in  \liv\cap L^2(0,T;W),\\
&\nonumber \s\in H^1(0,T;H)\cap C^0([0,T];V)\cap L^2(0,T;W)\,.
\end{align}
\noindent
(ii)\,\,\,{\em There is some constant $K_1^*>0$, which depends only on $R$ and the data
of the system,
such that for every $u\in\ur$ the associated strong solution $(\vp,\mu,\s)$ to} \eqref{ss1}--\eqref{ss5} {\em satisfies}
\beqa \label{ssbounds1}
&&\|\vp\|_{H^1(0,T;V)\cap L^\infty(0,T;H^3(\oma))}\,+\,  \|\Delta\vp\|_{L^2(0,T;\hzwei)}\,+\,
\|\mu\|_{\liv\cap L^2(0,T;H^2(\Omega))}\\[1mm]
&&\nonumber +\,
\|\s\|_{H^1(0,T;H)\cap C^0([0,T];V)\cap L^2(0,T;\hzwei)}\,\le\,K_1^*\,.
\eeqa
\noindent
(iii)\,\,{\em There is some constant $K_2^*>0$, which depends only on $R$ and the data of the problem, such that the
following holds true: whenever $u_i\in\ur$, $i=1,2$, are given and $(\vp_i,\mu_i,\s_i)$, $i=1,2$, are 
the associated solutions to the state system, then we have, for every $t\in [0,T]$,}
\beqa \label{stabu1}
&&\|\vp_1(t)-\vp_2(t)\|_{V'}\,+\,\|\vp_1-\vp_2\|_{\lzvt}\,+\,\|\s_1(t)-\s_2(t)\|_{V'}\\[1mm]
&&\nonumber +\,\|\s_1-\s_2\|_{\lzht}\,\le\,K_2^*\,\|u_1-u_2\|_{L^2(0,t;H)}\,.
\eeqa                                       

\vspace{3mm} \noindent
{\sc Proof:} \,\,\,In the following, we denote by $C_i>0$, $i\in\NN$, constants that depend only on $R$ and the data entering
the state system.  In \cite[Thms.\,1-3]{FGR}, it has been shown that the variational problem
\beqa \label{weak1}
&&\ioma\partial_t\vp(t)\,v\dx+\ioma\nabla\mu(t)\cdot\nabla v\dx\,=\,
\ioma P(\vp(t))(\s(t)-\mu(t))\,v\dx,\\[1mm]
\label{pier1}
&&\pier{\ioma\mu(t)\,v\dx\, =\, \ioma\nabla\vp(t)\cdot\nabla v\dx
+ \ioma F'(\vp(t))\,v\dx,}\\[1mm]
\label{weak2}
&&\ioma \partial_t\s(t)\,v\dx+\ioma\nabla\s(t)\cdot\nabla v\dx\,=\,
\ioma (u(t)-P(\vp(t))(\s(t)-\mu(t)))\,v\dx\,,
\eeqa
for all $v\in V$ and almost every $t\in (0,T)$, has in the homogeneous case $u\equiv 0$ a unique 
solution triple $(\vp,\mu,\s)$ which satisfies the initial conditions \eqref{ss5} and
has the regularity properties 
\beqa \label{regu1}
&&\vp\in \pier{ H^1(0,T;V)\cap{}} L^\infty(0,T;H^3(\oma)), \quad \pier{\mu} \in \liv,\\[1mm]
&&\s\in H^1(0,T;H)\cap \liv.\nonumber
\eeqa
A closer inspection of the proofs of \cite[Thms.\,1-3]{FGR} (in particular, 
\juerg{Eq.} \pier{ \eqref{weak2} is} nowhere 
differentiated with respect to time) reveals that only straightforward
modifications are needed to show that the system \pier{\eqref{weak1}--\eqref{weak2},} \eqref{ss5} 
has for every $u\in\ur$ a unique solution triple $(\vp,\mu,\s)$ 
which has the regularity properties \eqref{regu1} and satisfies, with some $C_1>0$, 
\beqa \label{2.14}
&&\|\vp\|_{H^1(0,T;V)\cap L^\infty(0,T;H^3(\oma))}\,+\, \|\mu\|_{\liv}
\,+\,
\|\s\|_{H^1(0,T;H)\cap L^\infty(0,T;V)}\,\le\,C_1\,.
\eeqa
Now observe that \eqref{2.14} \pier{and ({\bf H4}) imply} that $u-P(\vp)(\s-\mu)$ is bounded in $\lzq$.
Since $\sigma_0\in V$, parabolic regularity theory, applied to \eqref{weak2}, yields 
$\s\in L^2(0,T;W)$ and the bound for $\s$ stated in
\eqref{ssbounds1}. Moreover, Eq. \eqref{weak1} is for almost every $t\in (0,T)$ the
weak form of the elliptic problem 
$$-\Delta\mu(t)=P(\vp(t))(\s(t)-\mu(t))-\partial_t\vp(t) \quad\mbox{in }\,\oma,
\quad \dn\mu(t)=0\quad\mbox{on }\,\doma.$$
By \betti{({\bf H4}) and \eqref{2.14}, we have $\,P(\vp)(\sigma-\mu)-\partial_t\vp\in 
L^2(0,T;H)$, whence} we infer from elliptic regularity theory that $\mu(t)\in W$ for
almost every $t\in (0,T)$, as well as \bettinew{$\|\mu\|_{\lzhz}$ is uniformly bounded}. Moreover, we
have $\dn(\Delta\vp)=-\dn\mu+F''(\vp)\dn\vp=0$ almost everywhere on $\Sigma$, as
well as 
\beq
\label{Delta2}
\Delta^2\vp=-\Delta\mu + F''(\vp)\,\Delta\vp+F'''(\vp)|\nabla\vp|^2,
\eeq 
where\juerg{, due to \eqref{2.14} and ({\bf H5}),} the right-hand side 
is bounded in $L^2(Q)$ . We thus have
$\Delta\vp \in L^2(0,T;W)$ with \bettinew{bounded norm}, 
and the assertions (i) and (ii)
are proved\pier{.
The last} assertion (iii) can be shown in exactly the same way as the
stability result shown in the
proof of \cite[Thm.\,2]{FGR}. We therefore may skip the proof. \qed   

\vspace{5mm}\noindent
{\sc Remark 1}\quad Observe that standard embedding results 
(\pier{cf.~\cite[Sec.\,8, Cor.\,4]{simon}})
imply that $H^1(0,T;V)\cap L^\infty(0,T;\hzwei)$ is continuously embedded in
$C^0([0,T];H^s(\oma))$ for $0<s<2$. Consequently, $\vp\in C^0(\bq)$. Moreover, 
also owing to the continuity of the embeddings $V\subset L^6(\oma)$ and $\hzwei\subset\lio$ in
three dimensions of space, and invoking the general hypotheses ({\bf H4}) and ({\bf H5}), we may without loss of generality \betti{state} (possibly choosing a larger $K_1^*$) \betti{the following bounds:}
\beqa \label{ssbounds2}
&&\|\vp\|_{C^0(\bq)}\,+\,\|\nabla\vp\|_{L^\infty(Q)} \,+\,\max_{0\le i\le 4}\,
\|F^{(i)}(\vp)\|_{C^0(\bq)}\,+\,\max_{0\le i\le 2}\,\|P^{(i)}(\vp)\|_{\liq}\\[1mm]
&&\nonumber +\,\|\Delta\vp\|_{L^\infty(0,\juerg{T};L^6(\oma)\betti{)}\cap L^2(0,T;\lio)}
\,+\,\pier{\|\mu\|}_{L^\infty(0,T;L^6(\oma))
\cap L^2(0,T;\lio)}
\\[1mm]
&&\nonumber +\,\pier{\|\s\|}_{L^\infty(0,T;L^6(\oma))\cap L^2(0,T;\lio)}\,\le\, K^*_1,
\eeqa
for any solution $(\vp,\mu,\s)$ corresponding to some $u\in\ur$.

\vspace{2mm}\noindent
{\sc Remark 2}\quad  
A comparison argument yields that (by possibly choosing a larger $K_2^*$)  the estimate
\beq\label{stabu2}
\|\mu_1-\mu_2\|_{L^2(0,t;V')}\,\le\,K_2^*\,\|u_1-u_2\|_{L^2(0,t;H)}
\eeq
holds true for every $t\in [0,T]$ \betti{with $\mu_i:=-\Delta\vp_i
+F'(\vp_i)$, $i=1,2$}. In particular, the control-to-state operator $\cs$, $u\mapsto \cs(u):=(\vp,\mu,\s)$,
is well defined and Lipschitz continuous as a mapping from $\ur\subset\lzq$ into the space
$$(L^\infty(0,T,V')\cap \lzv)\times L^2(0,T;V')\times (L^\infty(0,T;V')\cap\lzh).$$ 

\vspace{2mm}
The stability results \eqref{stabu1}, \eqref{stabu2} are not sufficient for studying the control
problem $\PP$. We thus begin our analysis by proving stronger stability estimates. We have the following
result.
 
\vspace{5mm} \noindent
{\sc Theorem 2.2} \quad {\em Suppose that the hypotheses} $\as$ {\em are fulfilled. Then there exists a constant
$K_3^*>0$, which depends only on $R$ and the data of the system, such that the following holds true:
whenever $u_i\in\ur$, $i=1,2$, are given and $(\vp_i,\mu_i,\s_i)$, $i=1,2$, are 
the associated solutions to the state system \betti{\eqref{ss1}--\eqref{ss5}}, then we have\pier{, for every $t\in [0,T]$,}}
\beqa \label{stabu3}
&&\|\vp_1-\vp_2\|_{H^1(0,t;V')\cap L^\infty(0,t;V)\cap L^2(0,t;H^3(\oma))}\,+\,
\|\mu_1-\mu_2\|_{L^2(0,t;V)}\,\\[1mm]
&&\nonumber +\,\|\s_1-\s_2\|_{H^1(0,t;H)\cap C^0([0,t];V)\cap L^2(0,t;\hzwei)}\,
 \le\,K_3^*\,\|u_1-u_2\|_{\lzht}\,.
\eeqa          

\vspace{1mm} \noindent
{\sc Proof:} \,\,\,First, recall that the estimates \eqref{stabu1} and \eqref{stabu2} are valid.
We again denote by $C_i>0$, $i\in\NN$, constants that depend only on $R$ and the data entering
the state system. Putting $\vp:=\vp_1-\vp_2$, $\mu:=\mu_1-\mu_2$, $\s:=\s_1-\s_2$, and $u:=u_1-u_2$, we find \betti{the following} identities:
\beqa
\label{diff1}
&&\pt\vp-\Delta\mu\,=\,P(\vp_1)(\s-\mu)\,+\,(P(\vp_1)-P(\vp_2))(\s_2-\mu_2) \quad\mbox{a.\,e. in }\,Q,\\[2mm]
\label{diff2}
&&\mu\,=\,-\Delta\vp + F'(\vp_1)-F'(\vp_2)\quad\mbox{a.\,e. in }\,Q,\\[2mm]
\label{diff3}
&\quad&\pt\s-\Delta\s\,=\,u\,-\,P(\vp_1)(\s-\mu)\,-\,(P(\vp_1)-P(\vp_2))(\s_2-\mu_2)\quad\mbox{a.\,e. in }\,Q,\\[2mm]
\label{diff4}
&&\dn\vp=\dn\mu=\dn\s=0 \quad\mbox{a.\,e. on }\,\Sigma,\\[2mm]
\label{diff5}
&&\vp(0)=\s(0)=0 \quad\mbox{a.\,e. in }\,\oma.
\eeqa                         

We now establish the asserted stability estimate \eqref{stabu3} in a series of estimates in which we make
repeated use of the global bounds \eqref{ssbounds1} and \eqref{ssbounds2} without further reference. 

\vspace{2mm}
\noindent
\underline{\sc First estimate:} \,\,\,We multiply \eqref{diff1} by $\vp$ and \eqref{diff2} by $\Delta\vp$ and add the
resulting equations to obtain that
\begin{eqnarray*}
\vp\,\pt\vp -\vp\Delta\mu+ \mu\Delta\vp&\!\!=\!\!&P(\vp_1)(\s-\mu)\vp+(P(\vp_1)-P(\vp_2))(\s_2-\mu_2)\vp\\[1mm]
&&{}-|\Delta\vp|^2\,+\,(F'(\vp_1)-F'(\vp_2))\,\Delta\vp.
\end{eqnarray*}
Integration over $\oma\times [0,t]$ (where $t\in [0,T]$) and by parts \betti{yield}, 
using the mean value theorem, \betti{({\bf H4}), ({\bf H5})}, and Young's inequality, that
\beqa
\label{2.24}
&&\juerg{\frac 1 2\,\|\vp(t)\|_H^2\,+\,\frac 1 2 \txinto|\Delta\vp|^2\dx\ds}
\\[1mm]
&&\juerg{\,\le\, \betti{I_1}
\nonumber \,+\,\pier{C_2}\int_0^t(\|\s(s)\|_H^2\,+\,\|\vp(s)\|_H^2)\ds\,+\,\pier{C_3}
\,\betti{I_2},} 
\eeqa 
where, \pier{owing to ({\bf H4}), \eqref{ssbounds1}, \eqref{stabu1}, \eqref{ssbounds2}, \eqref{stabu2}, and   H\"older's inequality, we deduce~that
\begin{align*}
\betti{I_1} &:=\!\! \betti{\int_0^t\|\mu(s)\|_{V'} \|P(\vp_1(s))\vp(s)\|_V\ds }\\[1mm]
&\le \pier{\frac12}\int_0^t\|\mu(s)\|_{V'}^2\ds + \pier{\frac12}
\int_0^t\!\! \|P(\vp_1(s))\vp(s)\|_V^2\ds \\[1mm]
&\le\, \pier{C_{4}}\,\|u\|^2_{L^2(0,t;H)}\,+\, 
 \pier{\frac12}\txinto\!\! (|P(\vp_1)\vp|^2+|P(\vp_1)\nabla\vp+\vp \,P'(\vp_1)
\nabla\vp_1|^2)\dx\ds\\[1mm]
&\le\,\pier{C_{4}}\,\|u\|^2_{L^2(0,t;H)}\,+\, \pier{C_{5}}\int_0^t\|\vp(s)\|^2_V\ds\,\le\,\pier{C_{6}}\,\|u\|^2_{L^2(0,t;H)}\,,
\\[2mm]
\betti{I_2}&:=\!\!\txinto|\vp|^2\,(|\s_2|+|\mu_2|) \dx\ds\,\le\,
\int_0^t (\|\s_2(s)\|_6 +\|\mu_2(s)\|_6)\|\vp(s)\|^2_{12/5}\ds\\[1mm]
&\le\pier{C_7}\,\int_0^t\|\vp(s)\|_V^2\ds\,\le\,\pier{C_8}\,\|u\|^2_{L^2(0,t;H)}\,. 
\end{align*}%
}%
Combining the above inequalities, \pier{using again \eqref{stabu1}, and invoking well-known elliptic regularity results, we infer} the estimates
\beqa \label{2.26}
&&\|\vp\|_{L^\infty(0,t;H)\cap L^2(0,t;\hzwei)}\,\le \,\pier{C_{9}}\,\|u\|_{L^2(0,t;H)}
\quad\mbox{for all }\,t\in [0,T],\\[2mm] 
\label{2.27}
&&\|\mu\|_{\lzht}\,\le \,\pier{C_{10}}\,\|u\|_{L^2(0,t;H)}
\quad\mbox{for all }\,t\in [0,T],
\eeqa
where the second inequality follows \pier{from a} comparison in \eqref{diff2}, \betti{by \pier{applying}
once more the mean value theorem and ({\bf H5})}. 
 
\vspace{3mm}
\noindent
\underline{\sc Second estimate:} \,\,\,We now test \betti{\eqref{diff3}} by $\pt\s$ and use 
Young's inequality and the mean value theorem \betti{along with ({\bf H4})} to find that, for all 
$t\in [0,T]$,
\beqa\label{2.28}
&&\txinto |\pt\s|^2\dx\ds \,+\,\frac 12\,\|\nabla\s(t)\|_H^2\\[1mm]
&&\nonumber =\txinto \pt\s\,[u-P(\vp_1)(\s-\mu)-(P(\vp_1)-P(\vp_2))(\s_2-\mu_2)]\dx\ds\\[1mm]
&&\nonumber \le\,\pier{\frac 14}\txinto|\pt\s|^2\dx\ds\,+\,\pier{C_{11}}\txinto(u^2+\s^2+\mu^2)\dx\ds 
\,+\,\pier{I}, 
\eeqa
where, using H\"older's and Young's inequalities, as well as \eqref{ssbounds2}, 
\pier{we infer} \juerg{that}
\begin{eqnarray*}
I&\!\!:=\!\!& \txinto \pier{|\pt\s|\, |P(\vp_1)-P(\vp_2)|}\,(|\s_2|+|\mu_2|)\dx\ds\\[1mm]
&\!\!\le\!\!& \pier{C_{12}}\int_0^t\!\pier{\|\pt\s(s)\|_2\,\|\vp(s)\|_3}\,(\|\s_2(s)\|_6+\|\mu_2(s)\|_6)\,\ds\\[1mm]
&\!\!\le\!\!& \pier{\frac 14}\txinto |\pt\s|^2\dx\ds\,+\,\pier{C_{13}}\int_0^t\|\vp(s)\|_V^2\ds\,.
\end{eqnarray*}
In view of \betti{\eqref{stabu1}} and \eqref{2.27}, we thus obtain \betti{that} 
\beq\label{2.29}
\|\s\|_{H^1(0,t;H)\cap L^\infty(0,t;V)}\,\le\,\pier{C_{14}}\,\|u\|_{\lzht} \quad\mbox{for all }\,
t\in [0,T],
\eeq
whence, by comparison in \eqref{diff3}, \betti{and applying once more the mean value theorem together with ({\bf H4})}, also the \pier{bound}
\beq\label{2.30}
\|\s\|_{L^2(0,t;\hzwei)}\,\le\,\pier{C_{15}}\,\|u\|_{\lzht} \quad\mbox{for all }\,
t\in [0,T]
\eeq
\betti{follows} \pier{from the estimate of $\|\Delta\s\|_{L^2(0,t;H)}$ and the elliptic regularity theory.}

\vspace{3mm}
\noindent
\underline{\sc Third estimate:} \,\,\,Next, we insert $\mu$, given by \eqref{diff2}, in \eqref{diff1}
to find that 
\beqa\label{fourth}
\pt\vp+\Delta^2\vp&\!\!=\!\!&P(\vp_1)(\s-\mu)+(P(\vp_1)-P(\vp_2))(\s_2-\mu_2)+F''(\vp_1)
\Delta\vp\\[1mm]
&&\nonumber +\,(F''(\vp_1)-F''(\vp_2))\Delta\vp_2+F'''(\vp_1)(\nabla\vp_1+\nabla\vp_2)
\cdot\nabla\vp\\[1mm]
&&\nonumber +\,(F'''(\vp_1)-F'''(\vp_2))|\nabla\vp_2|^2\,.
\eeqa 
Testing this identity by $\,-\Delta\vp\in L^2(0,T;W)$, and using \betti{the integrations by parts,} Young's inequality, the
mean value theorem, \betti{({\bf H4}), ({\bf H5}),} 
and \eqref{ssbounds2}, we find that for any $t\in [0,T]$ we have 
\beqa\label{2.32}
&&\frac 12\,\|\nabla\vp(t)\|_H^2\,+\,\txinto|\nabla\Delta\vp|^2\dx\ds\\[1mm]
&&\nonumber \le\,\pier{C_{16}}\txinto|\vp|(|\s_2|+|\mu_2|+|\Delta\vp_2|)|\Delta\vp|\dx\ds\\[1mm]
&&\nonumber \quad +\,\pier{C_{17}}\txinto\Bigl(|\vp|^2+|\nabla\vp|^2+|\Delta\vp|^2+|\s|^2+|\mu|^2\Bigr)\dx\ds\,.
\eeqa
The first integral on the right-hand side, which we denote by $I$, can be estimated as follows:
\begin{eqnarray*}
I&\!\!\le\!\!&\int_0^t\|\vp(s)\|_3(\|\s_2(s)\|_6+\|\mu_2(s)\|_6+\|\Delta\vp_2(s)\|_6)\|\Delta\vp(s)\|_2\ds\\[1mm]
&\!\!\le\!\!&\pier{C_{18}}\int_0^t\Bigl(\|\vp(s)\|_V^2\,+\,\|\Delta\vp(s)\|_H^2\Bigr)\ds
\,\le\,\pier{C_{19}}\,\|u\|_{\lzht}^2,
\end{eqnarray*}
where the last inequality follows from \eqref{stabu1} and \eqref{2.26}. In conclusion, \betti{using once more \eqref{2.26} and \eqref{stabu1} in order to bound the second integral in \eqref{2.32},}  we have the estimate
\beq\label{2.33}
\|\vp\|_{L^\infty(0,t;V)\cap L^2(0,t;H^3(\oma))}\,\le\,\pier{C_{20}}\,\|u\|_{\lzht} \quad\mbox{for all }\,t\in [0,T].
\eeq
Comparison in \eqref{diff2} \betti{together with an application of the mean value theorem and of ({\bf H5})} then easily shows that also
\beq\label{2.34}
\|\mu\|_{\lzvt}\,\le\,\pier{C_{21}}\,\|u\|_{\lzht} \quad\mbox{for all }\,t\in [0,T]\,.
\eeq

\vspace{3mm}
\noindent
\underline{\sc Fourth estimate:} \,\,\,Finally, we test \eqref{diff1} by an
arbitrary function $v\in L^2(0,T;V)$ to obtain that
\beqa\nonumber
&&\left |\txinto \pt\vp\,v\dx\ds\right|\,\le\,\txinto|\nabla\mu||\nabla v|\dx\ds
\,+\,\pier{C_{22}}\txinto(|\s|+|\mu|)|v|\dx\ds\\[1mm]
&&\nonumber\quad+\,\pier{C_{23}}\int_0^t\|\vp(s)\|_4(\|\s_2(s)\|_4+\|\mu_2(s)\|_4)\,\|v(s)\|_2\dx\ds
\\[1mm]
&&\nonumber\le \,\pier{C_{24}}\left(\|\vp\|_{\lzvt}\,+\,\|\mu\|_{\lzvt}\,+\,
\|\s\|_{\lzht}\right)\,\|v\|_{\lzvt}\,.
\eeqa
In view of the above estimates, this implies that
\beq\label{2.35}
\|\pt\vp\|_{L^2(0,t;V')}\,\le\,\pier{C_{25}}\,\|u\|_{\lzht}\quad\mbox{for all }\,t\in [0,T].
\eeq
With \pier{\eqref{2.35}}, the assertion is completely proved.\qed

\section{Differentiability of the control-to-state operator}
\setcounter{equation}{0}

In this section, we establish a differentiability result for the control-to-state 
operator~$\cs$. To this end, we assume
that the  general assumptions $\as$ are satisfied. For arbitrary, but fixed $\bu\in\ur$, let
$(\bphi,\bmu,\bs)=\cs(\bu)$. We consider \pier{for} \betti{any} $h\in\lzq$ the linearized system
\beqa
\label{ls1}
&&\partial_t\xi-\Delta\eta\,=\,P'(\bphi)(\bs-\bmu)\,\xi\,+\,P(\bphi)(\rho-\eta) \quad\mbox{in }\,Q,\\[1mm]
\label{ls2}
&&\eta\,=\,-\Delta\xi+F''(\bphi)\,\xi \quad\mbox{in }\,Q,\\[1mm]
\label{ls3}
&&\partial_t\rho-\Delta\rho\,=\,-P'(\bphi)(\bs-\bmu)\,\xi\,-\,P(\bphi)(\rho-\eta)\,+\,h \quad\mbox{in }\,Q,\\[1mm]
\label{ls4}
&&\dn\xi \pier{{}= \dn \eta {}} =\dn\rho\,=\,0 \quad\mbox{on }\,\Sigma,\\[1mm]
\label{ls5}
&&\xi(0)=\rho(0)=0 \quad\mbox{in }\,\oma.
\eeqa
We expect the following: if the system \eqref{ls1}--\eqref{ls5} admits a unique solution 
$(\xi,\eta,\rho)=:(\xi^h,\eta^h,\rho^h)$ for every $h\in\lzq$ and the Fr\'echet derivative $D\cs(\bu)$ 
exists as a continuous linear mapping from $\ur$ into a 
suitable Banach space (which is yet to be determined), then we should have $D\cs(\bu)h=
(\xi^h,\eta^h,\rho^h)$. To this end, we first show the unique solvability of the linearized system
\eqref{ls1}--\eqref{ls5}. We have the following result.

\vspace{5mm}
\noindent
{\sc Theorem 3.1} \quad{\em The system} \eqref{ls1}--\eqref{ls5} {\em has a unique solution 
$(\xi,\eta,\rho)$ with}
\beqa
\label{reguls}
&&\xi\in H^1(0,T;V')\cap\betti{ L^\infty(0,T;V)\cap L^2(0,T;\pier{{}W\cap{}}H^3(\oma))}, \\[1mm]
&&\nonumber \pier{\eta\in L^2(0,T;V),}\quad \rho\in H^1(0,T;H)\cap C^0([0,T];V)\cap \betti{L^2(0,T;W)},
\eeqa
{\em which fulfills the conditions} \pier{\eqref{ls1}}--\eqref{ls5}  {\em almost everywhere 
in the respective sets, \pier{except for  \eqref{ls1} and the related boundary condition in 
\eqref{ls4} that are} fulfilled in the sense that, for almost every $t\in (0,T)$,}
\beqa
\label{ss1weak}
&&\langle \pt\xi(t),v\rangle_V\,+\,\ioma\nabla\eta(t)\cdot\nabla v\dx\,=\,
\ioma P(\bphi(t))(\rho(t)-\eta(t))\,v\dx\\[1mm]
&&\nonumber +\ioma P'(\bphi(t))(\bs(t)-\bmu(t))\,\xi(t)\,v\dx \quad\,\,\forall\,v\in V.
\eeqa   
{\em Moreover, there is some constant $K^*_4>0$, which depends only on $R$ and the data of the
state system, such that}
\beqa\label{conti}
&&\|\xi\|_{H^1(0,t;V')\cap L^\infty(0,t;V)\cap L^2(0,t;H^3(\oma))}
\,+\,\|\eta\|_{\lzvt}\\[1mm]
&&\nonumber
\,+\,\|\rho\|_{H^1(0,t;H)\pier{{}\cap C^0([0,t];V)}\cap L^2(0,t;\hzwei)}\, \le\,K_4^*\,\|h\|_{\lzht}\,\quad\mbox{for all }\,t\in [0,T].
\eeqa

\vspace{4mm}
\noindent
\pier{{\sc Remark 3} \quad Note that an equivalent formulation of \eqref{ss1weak}, which makes use of the abstract operator $A$ defined by \eqref{defA}, is 
\begin{align}
\label{ss1abs}
&\pt\xi(t) + A\eta(t)\,=\, P(\bphi(t))(\rho(t)-\eta(t))\\[1mm]
&\nonumber +P'(\bphi(t))(\bs(t)-\bmu(t))\,\xi(t) + \eta (t) \quad \hbox{in } V' ,  \ \hbox{ for a.e. } \, t\in (0,T).
\end{align}%
}

\noindent
{\sc Proof.} \quad We apply a Faedo-Galerkin approximation, using the family $\{w_j\}_{j\in\NN}\subset W$
of (appropriately ordered) eigenfunctions of the eigenvalue problem 
$$\,-\Delta w+w=\lambda w \quad\mbox{in $\oma$}, \quad \dn w=0 \quad\mbox{on $\doma$}, $$
as a Galerkin basis in $V$.
The family $\{w_j\}_{j\in\NN}$ forms a complete
orthonormal system in $(H, (\cdot,\cdot)_H)$ which is also orthogonal in 
$(V,(\cdot,\cdot)_V)$; moreover, we also have 
$(\Delta w_j,\Delta w_k)_H=0$ whenever $j\not =k$.

Let $n\in\NN$ be fixed. Putting ${\cal W}_n:={\rm span}\,\{w_1,...,w_n\}$, we then look for functions
of the form
\begin{equation*} 
\xi_n(x,t):=\sumk a_k^n(t)w_k(x), \quad \rho_n(x,t):=\sumk b_k^n(t)w_k(x), \quad
\eta_n(x,t):=\sumk c_k^n(t)w_k(x), 
\end{equation*}
which satisfy the following approximating problem:
\beqa
\label{gal1}
&&\ioma \pt\xi_n(t)v\dx \,+\ioma\nabla\eta_n(t)\cdot\nabla v\dx \,-\ioma P(\bphi(t))(\rho_n(t)-\eta_n(t))v\dx\\[1mm]
\nonumber
&&=\,\ioma P'(\bphi(t))(\bs(t)-\bmu(t))\xi_n(t)v\dx\quad\forall\,v\in {\cal W}_n,\\[2mm]
\label{gal2}
&&\ioma\eta_n(t)v\dx\,=\,\ioma \nabla\xi_n(t)\cdot\nabla v\dx\,+\ioma F''(\bphi(t))\xi_n(t)v\dx
\quad\forall\,v\in {\cal W}_n,\\[2mm]
\label{gal3}
&&\ioma \pt\rho_n(t)v\dx \,+\ioma\nabla\rho_n(t)\cdot\nabla v\dx \,+\ioma P(\bphi(t))(\rho_n(t)-\eta_n(t))v\dx\\[1mm]
\nonumber
&&=\ioma [-P'(\bphi(t))(\bs(t)-\bmu(t))\xi_n(t) \,+\,h(t)]v\dx\quad\forall\,v\in {\cal W}_n,\\[2mm]
\label{gal4}
&&\xi_n(0)=\rho_n(0)=0.
\eeqa

\noindent By insertion of $v=w_j$, $1\le j\le n$, in \eqref{gal2}, it is easily seen that 
the unknowns $c_j^n$ can be expressed explicitly in terms of the unknowns $a_1^n
,\ldots,a_n^n$; hence, by inserting $v=w_j$, $1\le j\le n$, in \eqref{gal1} and 
\eqref{gal3}, it turns out that 
the system \eqref{gal1}--\eqref{gal4}
is in fact equivalent to a Cauchy problem for a linear system of $2n$ 
first-order ordinary differential equations in the $2n$ 
unknowns $a_j^n,b_j^n$, $1\le j\le n$, in which, owing to ({\bf H4}) and 
({\bf H5}), and to the fact that 
$\bphi\in C^0(\bq)$ (recall \eqref{ssbounds2}), 
all of the occurring coefficient functions
belong to $L^\infty(0,T)$. By Carath\'eodory's theorem, this linear system has a unique
solution $(a_1^n,\ldots,a_n^n,b_1^n,\ldots,b_n^n)\in (W^{1,\infty}(0,T))^{2n}$, which 
specifies the unique solution $(\xi_n,\eta_n,\rho_n)\in (W^{1,\infty}(0,T;{\cal W}_n))^3$
to \eqref{gal1}--\eqref{gal4}.
    
We now aim to derive a number of a priori estimates for the approximations. To this end, we denote by $C_i$, $i\in\NN$, positive constants \betti{(possibly different from the ones used in the previous section)} that may depend on the data of the problem, but 
not on $n\in\NN$.  

\noindent \underline{\sc First estimate:}  \,\,\,Let $t\in [0,T]$ be arbitrary. We insert $v=\eta_n$
in \eqref{gal1}, $v=-\pt\xi_n$ in \eqref{gal2}, as well as $v=\rho_n$ in \eqref{gal3}, and add
the resulting identities. Integrating over $\oma\times [0,t]$ and by parts, we obtain that
\beqa\label{3.13}
&&\quad\frac 12\,(\|\rho_n(t)\|_H^2+ \|\nabla\xi_n(t)\|_H^2)\,+\txinto(|\nabla\rho_n|^2+|\nabla\eta_n|^2)\dx\ds\\[1mm]
&&\nonumber +\txinto P(\bphi)\pier{(\rho_n- \eta_n)^2}\dx\ds \,=\,-\txinto\pt\xi_n\,F''(\bphi)\,\xi_n\dx\ds\\[1mm]
&&\nonumber \pier{\,+
\txinto h\,\rho_n\dx\ds }\,+\txinto P'(\bphi)\,(\bs-\bmu)\,\xi_n\,(\eta_n-\rho_n)\dx\ds\,.
\eeqa 
Owing to ({\bf H4}), the last integral on the left is nonnegative. We denote the integrals on the right by
$I_i$, \pier{$i=1,2,3$}, in that order, and estimate them indivually. Clearly, by Young's inequality,
\pier{it turns out that} 
\beq
\label{3.14}
|I_2|\,\le\,\frac 12\, (\|h\|_{\lzht}^2\,+\,\|\rho_n\|^2_{\lzht})\,.
\eeq 
Moreover, by virtue of \eqref{ssbounds2} and H\"older's and Young's inequalities, \pier{we have that} 
\beqa\label{3.15}
&&|I_1|\,=\,\left|-\frac 12\ioma F''(\bphi(t))|\xi_n(t)|^2\dx\,+\,\frac 12 \txinto F'''(\bphi)\,\pt\bphi\,\xi_n^2\dx\ds\right|\\[1mm]
&&\nonumber \le \,\frac {\pier{K_1^*}}2\,\|\xi_n(t)\|_H^2\,+\,C_1\pier{\int_0^t}\|\pt\bphi(s)\|_6\,\|\xi_n(s)\|_3
\,\|\xi_n(s)\|_2\ds\\[1mm]
&&\nonumber \le\,\frac {\pier{K_1^*}}2\,\|\xi_n(t)\|_H^2\,+\,C_2\,\pier{\int_0^t}\Bigl(\|\xi_n(s)\|_V^2\,+\,
\|\pt\bphi(s)\|_V^2\,\|\xi_n(s)\|_H^2\Bigr)\ds\,,
\eeqa
where we notice that the mapping $s\mapsto \|\pt\bphi(s)\|_V^2$ belongs to $L^1(0,T)$, \betti{due to \eqref{ssbounds1}}. Next, we observe that, for every $\pier{\gamma}>0$ (to be chosen later),
\beqa\label{3.17}
|\pier{I_3}|&\!\!\le\!\!&\pier{C_3}\int_0^t(\|\bs(s)\|_6\, +\,\|\bmu(s)\|_6)\,\|\xi_n(s)\|_3\,(\|\eta_n(s)\|_2+
\|\rho_n(s)\|_2)\ds\\[1mm]
\nonumber &\!\!\le\!\!&\gamma\,(\|\eta_n\|_{\lzht}^2\,+\,\|\rho_n\|^2_{\lzht})\,+\,\frac{\pier{C_4}}\gamma
\pier{ \|\xi_n(s)\|_{\lzvt}^2}\,.
\eeqa

We still need estimates for the $L^2(Q)$ norm of $\eta_n$ and for the $L^\infty(0,T;H)$
norm of $\xi_n$. To obtain the former, we insert $v=\eta_n$ in \eqref{gal2}, integrate
over $[0,t]$ and by parts, and apply
Young's inequality and \eqref{ssbounds2} to deduce that
\beqa
\nonumber
&&\txinto\eta_n^2\dx\ds\,=\,\txinto\eta_n(-\Delta\xi_n+F''(\bphi)\xi_n)\dx\ds\\[1mm]
&&\nonumber \le\,\frac 12\,\txinto\eta_n^2\dx\ds\,+\betti{\,\frac 12\txinto\pier{(|\Delta\xi_n|
+\, \pier{K_1^*}\, |\xi_n|)^2}}\dx\ds,
\eeqa
whence
\beq\label{3.18}
\|\eta_n\|^2_{\lzht}\,\le\,\pier{2\txinto|\Delta\xi_n|^2\dx\ds\,+\,
C_5\,\|\xi_n\|^2_{\lzht}\,}.
\eeq
To derive the missing estimates, we finally insert $v=\pier{(2+K_1^*)}\,\xi_n$ 
in \eqref{gal1} and $v=\pier{(2+K_1^*)}\Delta\xi_n$
in \eqref{gal2} and add the resulting equations. Integration over $[0,t]$ and by parts then yields that 
\betti{%
\beqa\label{3.19}
&&\pier{\left(1+ \frac{K_1^*}2 \right)}\,\|\xi_n(t)\|_H^2\,+\,\pier{(2+K_1^*)}\txinto|\Delta\xi_n|^2\dx\ds
\\[1mm]
&&\nonumber
\,=\, \txinto \pier{(2+K_1^*)}\, P(\bphi)\,(\rho_n-\eta_n)\, \xi_n\dx\ds\\[1mm]
&&\nonumber\quad\,\,\,+\,\txinto \pier{(2+K_1^*)} \,P'(\bphi)(\bs-\bmu)\xi_n^2\dx\ds\,
\\[1mm]
&&\nonumber\quad\,\,\,+\,\txinto \pier{(2+K_1^*)}\,
F''(\bphi)\,\xi_n\,\Delta\xi_n\dx\ds.
\eeqa
}%
We denote the integrals on the right-hand side by \pier{$I_4,\, I_5,\, I_6$}, in this 
order, and estimate them individually. First, we obviously have
\beq \label{3.20}
|\pier{I_6}|\,\le\,\pier{\frac{K_1^*}2 }\txinto|\Delta\xi_n|^2\dx\ds\,+\,\pier{C_6}\txinto|\xi_n|^2\dx\ds\,.
\eeq 
Moreover, owing to H\"older's and Young's inequalities and \eqref{ssbounds2}, \pier{we infer that}
\beqa\label{3.21}
\quad|\pier{I_5}|&\!\!\le\!\!&C_7\int_0^t(\|\bs(s)\|_6+\|\bmu(s)\|_6)\,\|\xi_n(s)\|_3\,\|\xi_n(s)\|_2\ds
\,\le\,C_8\,\|\xi_n\|^2_{\lzvt},
\eeqa
and, using Young's inequality once more,
\beq\label{3.22}
|\pier{I_4}|\,\le\,\gamma\,(\|\eta_n\|^2_{\lzht}\,+\,\|\rho_n\|^2_{\lzht})\,+\,\frac{C_9}\gamma
\,\|\xi_n\|^2_{\lzht}\,.
\eeq

\noindent
Now, we \pier{take the sum of \eqref{3.13}, \eqref{3.18}, and \eqref{3.19}\juerg{. Then, on account of \eqref{3.14}--\eqref{3.17} and  of} \eqref{3.20}--\eqref{3.22},  we can choose $\gamma $ small enough 
(namely, $\gamma <1/2$)  and apply Gronwall's lemma in order to} find the estimate
\beqa\label{3.23}
&&\|\xi_n\|_{L^\infty(0,t;V)\cap L^2(0,t;\hzwei)}\,+\,\|\eta_n\|_{\lzvt}
\,+\,\|\rho_n\|_{L^\infty(0,t;H)\cap\lzvt}\\[1mm]
&&\nonumber
\le \,C_{10}\,\|h\|_{\lzht}\,\quad\mbox{for all }\,t\in[0,T].
\eeqa

\vspace{4mm}
\noindent
\underline{\sc Second estimate:}  \,\,\,Let $t\in [0,T]$ be arbitrary. Observing that
$v=\Delta^2\xi_n\in {\cal W}_n$, we obtain from \eqref{gal2}, using integration by parts, 
Young's inequality, and \eqref{ssbounds2}, that
\begin{eqnarray*}
&&\txinto |\nabla\Delta\xi_n|^2\dx\ds=\txinto\nabla\xi_n\cdot\nabla\Delta^2\xi_n\dx\ds
=\txinto\Delta^2\xi_n(\eta_n-F''(\bphi)\xi_n)\dx\ds\\[1mm]
&&=\,-\txinto\nabla\Delta\xi_n\cdot(\nabla\eta_n-\xi_n\,F'''(\bphi)\nabla\bphi-F''(\bphi)
\nabla\xi_n)\dx\ds\\[1mm]
&&\le\,\frac 12\txinto|\nabla\Delta\xi_n|^2\dx\ds\,+\,C_{11}\txinto(|\nabla\eta_n|^2
+|\xi_n|^2+|\nabla\xi_n|^2)\dx\ds\,,
\end{eqnarray*}  
and it follows from \eqref{3.23} that
\beq\label{3.24}
\|\xi_n\|_{L^2(0,t;H^3(\oma))}\,\le\,C_{12}\,\|h\|_{\lzht} \quad\mbox{for all }\,
t\in [0,T].
\eeq

\vspace{4mm}
\noindent
\underline{\sc Third estimate:}  \,\,\,Let $t\in [0,T]$ be arbitrary.
We insert $v=\pt\rho_n$ in \eqref{gal3} and integrate over $[0,t]$ and by parts.
Using \eqref{ssbounds2} and Young's and H\"older's inequalities, we obtain 
for every $\gamma>0$ the 
estimate
\begin{eqnarray*}
&&\txinto|\pt\rho_n|^2\dx\ds\,+\,\frac 12\,\|\nabla\rho_n(t)\|_H^2
\,\le\,\,\betti{\gamma}\txinto|\pt\rho_n|^2\dx\ds\\[1mm]
&&+\,\frac{C_{13}}\gamma\txinto(\rho_n^2+\eta_n^2+h^2)\dx\ds\,+\,C_{14}\,I
\end{eqnarray*}
\pier{holds,} where
\begin{eqnarray*}
I&\!\!:=\!\!&\int_0^t(\|\bs(s)\|_6+\|\bmu(s)\|_6)\,\|\xi_n(s)\|_3\,
\|\pt\rho_n(s)\|_2\ds\\[1mm]
&\!\!\le\!\!&\gamma\txinto|\pt\rho_n|^2\dx\ds\,+\,\frac{C_{15}}\gamma\,
\|\xi_n\|_{\lzvt}^2\,,
\end{eqnarray*}
Choosing \betti{$1/2>\gamma>0$}, we obtain from \eqref{3.23} that
\beq\label{3.25}
\|\rho_n\|_{H^1(0,t;H)\cap L^\infty(0,t;V)}\,\le\,C_{16}\,\|h\|_{\lzht}
\quad\mbox{for all }\,t\in [0,T].
\eeq
Similar reasoning, using $v=-\Delta\rho_n$ in \eqref{gal3}, yields that also
\begin{equation*}
\txinto|\Delta\rho_n|^2\dx\ds\,\le\,C_{17}\,\|h\|_{\lzht}^2\,,
\end{equation*}
so that
\beq\label{3.26}
\|\rho_n\|_{L^2(0,t;H^2(\oma))} \,\le\,C_{18}\,\|h\|_{\lzht}\,.
\eeq
In conclusion, we have shown the estimate
\beqa\label{3.27}
&&
\|\xi_n\|_{L^\infty(0,t;V)\cap L^2(0,t;H^3(\Omega))}\,+\,\|\eta_n\|_{\lzvt}\\[1mm]
&&\nonumber\,+\,\|\rho_n\|_{H^1(0,t;H)\betti{{}\cap {}\pier{{}C^0([0,t];V)}}
\cap L^2(0,t;H^2(\oma))}
\,\le\,C_{19}\,\|h\|_{\lzht}\,.
\eeqa

\vspace{4mm}
\noindent
\underline{\sc Conclusion of the proof:}  \,\,\,It follows from \eqref{3.27} 
that there are functions $(\xi,\eta,\rho)$
such that, possibly only for a subsequence which is again indexed by $n$,
\begin{eqnarray*}
\xi_n\to\xi&&\mbox{weakly in } \,L^2(0,T;H^3(\oma))\\
&&\mbox{and weakly star in }\,L^\infty(0,T;V),\\[1mm]
\eta_n\to\eta&&\mbox{weakly in }\,L^2(0,T;V),\\[1mm]
\rho_n\to\rho&&\mbox{weakly in }\,H^1(0,T;H)\cap L^2(0,T;\hzwei).
\end{eqnarray*}
From the semicontinuity properties of the involved norms, we can infer
that the estimate \eqref{3.27} holds true for $(\xi,\eta,\rho)$ in place
of $(\xi_n,\eta_n,\rho_n)$, and it is easily seen that $(\xi,\eta,\rho)$ satisfies
\eqref{ls2} and \eqref{ls3} almost everywhere in $Q$. Moreover, we have $\rho(0)=0$
almost everywhere in $\oma$ and $\dn\xi=\dn\rho=0$ almost everywhere on $\Sigma$,
and it follows that, for every $v\in H^1(0,T;V)$ with $v(T)=0$, it holds the
identity
\begin{align*}
&-\int_0^T\!\!\!\ioma\xi\,\pt v\dx\dt\,+\int_0^T\!\!\!\ioma \nabla\eta\cdot\nabla v\dx\dt
\\[1mm]
&\quad =\,\int_0^T\!\!\!\ioma \left( P'(\bphi)(\bs-\bmu)\,\xi\,+\,P(\bphi)
(\rho-\eta)\right)v\dx\dt\,,
\end{align*}
which implies that \eqref{ss1weak} and $\xi(0)=0$ hold true. 
\pier{Indeed, we also recover that} $\xi\in H^1(0,T;V')$, and \pier{in addition comparison yields} 
$$\|\xi\|_{H^1(0,t;V')}\,\le\,C_{20}\,\|h\|_{\lzht}\,,$$
so that \eqref{conti} is shown. 

To prove
uniqueness, we write the system \eqref{ss1weak}, \eqref{ls2}--\eqref{ls5} for two
solutions $(\xi_i,\eta_i,\rho_i)$, $i=1,2$, and subtract the equations. Then
$\xi:=\xi_1-\xi_2$, $\eta:=\eta_1-\eta_2$, $\rho:=\rho_1-\rho_2$\, satisfy the system
\eqref{ss1weak}, \eqref{ls2}--\eqref{ls5} with $h\equiv 0$. Now notice that, up to obvious 
modifications which are necessary due to the fact that
we only have $\pt\xi\in L^2(0,T;V')$, the estimates leading to \eqref{3.23} can
be repeated. \pier{We point out, in particular, that all the three terms of the equation 
\eqref{ls2} belong to $ L^2(0,T;V)$ thanks to \eqref{reguls}, ({\bf H5}) and \eqref{reguss}. Then, since}
$h\equiv 0$ in this case, we must have $\xi=\eta=\rho=0$ and thus uniqueness.\qed 

We are now in a position to establish the Fr\'echet differentiability of
the control-to-state operator. We have the following result.

\vspace{3mm}\noindent
{\sc Theorem 3.2} \quad{\em Suppose that the assumptions $\as$ are satisfied. Then
the control-to-state mapping ${\cal S}$ is Fr\'echet differentiable in $\ur$ as
a mapping from $\lzq$ into the space 
\beqa \label{dspace}
{\cal Y}&\!\!:=\!\!&\Bigl(H^1(0,T;W')\cap L^\infty(0,T;H)\cap L^2(0,T;W)\Bigr)
\times L^2(Q)\\[1mm]
&&\nonumber\quad\times\Bigl(H^1(0,T;H)\cap L^2(0,T;\hzwei)\Bigr).
\eeqa
 Moreover, for any $\bu\in\ur$
the Fr\'echet derivative $D{\cal S}(\bu)\in {\cal L}(\lzq,{\cal Y})$ is defined as 
follows: for any $h\in\lzq$, we have $D{\cal S}(\bu)h=(\xi^h,\eta^h,\rho^h)$, where
$(\xi^h,\eta^h,\rho^h)$ is the unique solution to the linearized system} 
\eqref{ls1}--\eqref{ls5} {\em associated with $h$}.  

\vspace{2mm}\noindent
{\sc Proof:} 
Let $\bu\in\ur$ be arbitrary, and $(\bphi,\bmu,\bs)= \cs(\bu)$. 
Since $\ur$ is open, there is some $\Lambda>0$ such that $\bu+h\in\ur$ 
whenever $h\in\lzq$ and $\|h\|_{\lzq}\,\le\,\Lambda$. In the following, we
only consider such variations $h\in\lzq$. We put $(\vp^h,\mu^h,\s^h):=\cs(\bu+h)$,
and we denote by $(\xi^h,\eta^h,\rho^h)$ the unique solution to the linearized
system \eqref{ls1}--\eqref{ls5} associated with $h$. Notice that by \eqref{conti} the 
linear mapping $h\mapsto
(\xi^h,\eta^h,\rho^h)$ is continuous between the spaces $\lzq$ and ${\cal Y}$.
We now define
$$\psi^h:=\vp^h-\bphi-\xi^h,\quad \zeta^h:=\mu^h-\bmu-\eta^h,\quad \chi^h:=
\s^h-\bs-\rho^h.$$
According to Theorem 2.1 and Theorem 3.1, we have the regularity
\beqa\label{3.29}
&&\psi^h\in H^1(0,T;V')\cap L^\infty(0,T;V)\cap L^2(0,T;H^3(\oma)), \quad
\zeta^h \in L^2(0,T;V),\\[1mm]
&&\nonumber \chi^h\in H^1(0,T;H)\cap C^0([0,T];V)\cap L^2(0,T;\hzwei)\,.
\eeqa
Note also that $(\vp^h,\mu^h,\s^h)$ and $(\bphi,\bmu,\bs)$ satisfy the global
bounds \eqref{ssbounds1} and \eqref{ssbounds2}. \pier{Let us point out that $\psi^h$ is at least strongly continuous from $[0,T]$ to $H$ (see, e.g., \cite[Sec.\,8, Cor.\,4]{simon}).}

According to the definition of Fr\'echet differentiability, it suffices to show 
that there exists an increasing function $Z:(0,+\infty)\to (0,+\infty)$ with
$\,\lim_{\lambda\searrow 0}\,Z(\lambda)/\lambda^2=0\,$ and 
\beq\label{3.30}
\|(\psi^h,\zeta^h,\chi^h)\|_{\cal Y}^2\,\le\,Z(\|h\|_{L^2(Q)})\,.
\eeq
Now observe that $(\psi^h,\zeta^h,\chi^h)$ is a solution to the following
problem:
\beqa
\label{fs1}
&&\quad\langle \pt\psi^h(t),v\rangle_V \,+\,\ioma\nabla\zeta^h(t)\cdot\nabla v\dx
\,=\,\ioma \big[P(\vp^h)(\s^h-\mu^h)-P(\bphi)(\bs-\bmu)\\[2mm]
&&\nonumber -\,P'(\bphi)(\bs-\bmu)\,\xi^h - P(\bphi)(\rho^h-\eta^h)\big]
(t)\,v\dx \quad\mbox{for all $v\in V$ and a.\,e. $t\in (0,T)$,}\\[2mm]
\label{fs2}
&&\quad\zeta^h\,=\,-\Delta\psi^h+F'(\vp^h)-F'(\bphi)-F''(\bphi)\,\xi^h \quad\mbox{a.\,e
in }\,Q,\\[2mm]
\label{fs3}
&&\quad\pt\chi^h-\Delta\chi^h\,=\,-P(\vp^h)(\s^h-\mu^h)\,+\,P(\bphi)(\bs-\bmu)\,+\,
P'(\bphi)(\bs-\bmu)\,\xi^h\\[1mm]
&&\nonumber\hspace*{3.3cm} + P(\bphi)(\rho^h-\eta^h)\quad\mbox{a.\,e. in }\,Q,\\[2mm]
\label{fs4}
&&\quad\dn\psi^h=\dn\chi^h=0 \quad\mbox{a.\,e. on }\,\Sigma,
\\[2mm]
\label{fs5}
&&\quad \psi^h(0)=\chi^h(0)=0 \quad\mbox{a.\,e. in }\,\oma.
\eeqa
We also note that
a straightforward computation, using Taylor's theorem with integral remainder, yields the
identities
\beqa\label{taylor1}
&&F'(\vp^h)-F'(\bphi)-F''(\bphi)\,\xi^h\,=\,F''(\bphi)\,\psi^h\,+\,R_1^h\,(\vp^h-\bphi)^2,\\[1mm]
\label{taylor2}
&&P(\vp^h)(\s^h-\mu^h)-P(\bphi)(\bs-\bmu)-P'(\bphi)(\bs-\bmu)\,\xi^h - P(\bphi)(\rho^h-\eta^h)\\[1mm]
&&\nonumber =\,P(\vp^h)\,\s^h-P(\vp^h)\,\mu^h-P(\bphi)\,\bs+ P(\bphi)\,\bmu
-P'(\bphi)\,\bs\,\xi^h+P'(\bphi)\,\bmu\,\xi^h\\[1mm]
&&\nonumber\quad -P(\bphi)(\s^h-\bs-\chi^h)+P(\bphi)(\mu^h-\bmu-\zeta^h)\\[1mm]
&&\nonumber =P(\bphi)\,\chi^h+(P(\vp^h)-P(\bphi))(\s^h-\bs)+ (P(\vp^h)-P(\bphi)
-P'(\bphi)\,\xi^h)\,\bs\\[1mm]
&&\nonumber\quad -P(\bphi)\,\zeta^h-(P(\vp^h)-P(\bphi))(\mu^h-\bmu)- (P(\vp^h)-P(\bphi)
-P'(\bphi)\,\xi^h)\,\bmu\\[1mm]
&&\nonumber =\, P(\bphi)(\chi^h-\zeta^h)\,+\,(P(\vp^h)-P(\bphi))\,[(\s^h-\bs)-(\mu^h-\bmu)]\\[1mm]
&&\nonumber\quad +\,P'(\bphi)\,(\bs-\bmu)\,\psi^h\,+\,(\bs-\bmu)\,R_2^h\,(\vp^h-\bphi)^2\,=:\,Q^h\,,
\eeqa
where 
\beq\label{taylor3}
R_1^h\,=\int_0^1(1-z)\,F'''(\bphi+z(\vp^h-\bphi))\dz,\quad 
R_2^h\,=\int_0^1(1-z)\,P''(\bphi+z(\vp^h-\bphi))\dz\,.
\eeq

In the following estimates, we denote by $C_i$, $i\in\NN$, positive constants \betti{(possibly different from the ones used in the previous sections)} which may depend
on the data of the state system but not on \juerg{$h\in\lzq$} with $\|h\|_{\pier{L^2(Q)}} \,\le\,\Lambda$. 
For the sake of a better readability, we will often suppress the superscript $h$ in the functions
$(\psi^h,\zeta^h,\chi^h)$ during the estimates and only write them in the final
estimate in each step. 
We first notice
that, thanks to ({\bf H4}), ({\bf H5}), and \eqref{ssbounds2},
\beqa
\label{3.39}
\left\|R_1^h\right\|_{\liq} \,+\,\left\|\pier{R_2^h}\right\|_{\liq}\,\le\,C_1.
\eeqa
We also recall \betti{that} the \pier{inequalities} \eqref{embed1}--\eqref{GN1}, the global bounds \eqref{ssbounds1}, \eqref{ssbounds2}, the global stability
estimates \eqref{stabu1}, \eqref{stabu2}, \eqref{stabu3}, and the properties \eqref{ruleA1},
\eqref{ruleA2} satisfied by the Riesz isomorphism $A$ introduced in \eqref{defA} \betti{will be frequently used in the sequel without mentioning them}. We begin our analysis
by proving some preparatory $L^2$ estimates. We have, for every $t\in [0,T]$:

\begin{align}
\label{neu1}
&\int_0^t\|R_1^h(s)\,(\vp^h(s)-\bphi(s))^2\|_H^2\ds\,\le\,C_2\txinto|\vp^h-\bphi|^4\dx\ds\\[1mm]
&\nonumber\le\,C_2\,\int_0^t\|\vp^h(s)-\bphi(s)\|_\infty^2\,\|\vp^h(s)-\bphi(s)\|_H^2\ds
\,\le\,C_3\,\|h\|^4_{\lzht}\,.
\end{align}
Moreover, owing to \eqref{embed2}, \eqref{ssbounds1}, \eqref{ssbounds2}, ({\bf H4}), and
\eqref{stabu3}, \pier{we infer that}
\begin{align}
\label{neu2}
&\int_0^t\|(P(\vp^h(s))-P(\bphi(s)))((\s^h(s)-\bs(s))-(\mu^h(s)-\bmu(s)))\|_H^2\ds\\[1mm]
&\nonumber\le\,C_4\int_0^t\|P(\vp^h(s))-P(\bphi(s))\|_V^2\,\left(\|\s^h(s)-\bs(s)\|_V^2
+\|\mu^h(s)-\bmu(s)\|_V^2\right)\ds\nonumber\\[1mm]
&\nonumber\le\,C_5\int_0^t\|\vp^h(s)-\bphi(s)\|_V^2\,\left(\|\s^h(s)-\bs(s)\|_V^2
+\|\mu^h(s)-\bmu(s)\|_V^2\right)\ds\\[1mm]
&\nonumber\le\,C_6\,\|h\|^4_{\lzht}\,,
\end{align}
as well as, using H\"older's inequality, \eqref{ssbounds2}, and \eqref{stabu3},
\begin{align}
\label{neu3}
&\int_0^t\|(\bs(s)-\bmu(s))\,R^h_2(s)\,(\vp^h(s)-\bphi(s))^2\|_H^2\ds\\[1mm]
&\nonumber\le \,C_7\txinto\left(|\bs|^2+|\bmu|^2\right)|\vp^h-\bphi|^4\dx\ds\\[1mm]
&\nonumber\le\,C_8\int_0^t\left(\|\bs(s)\|_6^2+\|\bmu(s)\|^2_6\right)\,
\|\vp^h(s)-\bphi(s)\|_6^4\ds\,\le\,C_9\,\|h\|^4_{\lzht}\,.
\end{align}
Moreover, we have that
\begin{align}
\label{neu4}
&\int_0^t\|P'(\bphi(s))(\bs(s)-\bmu(s))\psi(s)\|_H^2\ds
\,\le\,C_{10}\txinto\betti{\left(|\bs|^2+|\bmu|^2\right)|\psi|^2}\dx\ds\\[1mm]
&\nonumber\le\,C_{11}\int_0^t\!\!\left(\|\bs(s)\|_6^2+\|\bmu(s)\|_6^2\right)
\|\psi(s)\|_3^2\ds
\,\le\,C_{12}\int_0^t\!\!\|\psi(s)\|_V^2\ds\,.
\end{align}

\vspace{4mm}
\noindent
\underline{\sc First estimate:}  \,\,\,First, we observe that Eqs. \eqref{fs1}--\eqref{fs3} can be rewritten in the form 
\beqa\label{3.44}
&&\pt\psi + A\zeta = Q^h+\zeta, \quad \zeta=A\psi + F''(\bphi)\psi
+ R_1^h\,(\vp^h-\bphi)^2-\psi,\\[1mm]
&&\nonumber \pt\chi+A\chi=-Q^h+\chi \quad \pier{\hbox{in } V', \ \hbox{ a.e. in } (0,T), }
\eeqa 
\betti{where $Q^h$ is defined in \eqref{taylor2}}. 

We now test the first equation in \eqref{3.44} by $A^{-1}\psi$, the third by $A^{-1}\chi$ and add the
resulting identities. Using \eqref{ruleA1} and \eqref{ruleA2}, we easily deduce that, for any $t\in [0,T]$,
\beqa\label{3.45}
&&\frac 12\left(\|\psi(t)\|_{V'}^2\,+\,\|\chi(t)\|_{V'}^2\right)+\txinto\left(|\nabla\psi|^2+|\chi|^2\right)
\dx\ds\\[1mm]
&&\nonumber\le\,-\txinto F''(\bphi)\psi^2\dx\ds\,-\txinto R_1^h(\vp^h-\bphi)^2\,\psi\dx\ds
\,+\int_0^t\|\chi(s)\|_{V'}^2\ds\\[1mm]
&&\nonumber\quad \,\,+\int_0^t\langle Q^h(s),A^{-1}\psi(s)-A^{-1}\chi(s)\rangle_V\ds
\,+\int_0^t\langle\zeta(s),A^{-1}\psi(s)\rangle_V\ds\,.
\eeqa

\noindent
We denote the first, second, fourth, and fifth integral on the right-hand side 
by $I_1,I_2,I_3,I_4$, in this
order. Using \eqref{ssbounds2}, \eqref{neu1}, and Young's inequality, we have
\beq\label{3.46}
|I_1|\,+\,|I_2|\,\le\,C_{13}\left(\|\psi\|^2_{\lzht}\,+\,\|h\|^4_{\lzht}\right).
\eeq 

\noindent
Using \eqref{3.44} and \eqref{neu1}, we also have
\beqa\label{3.47}
|I_4|&\!\!=\!\!&\left|\int_0^t\langle (A\psi-\psi+F''(\bphi)\psi+R_1^h(\vp-\bphi)^2)(s),A^{-1}\psi(s)\rangle_V\ds
\right|\\[1mm]
\nonumber&\!\!\le\!\!&C_{14}\int_0^t\left(\|\psi(s)\|_H^2+\|\psi(s)\|_{V'}^2 \right)\ds 
\,+\,C_{15}\,\|h\|^4_{\lzht}\,.
\eeqa

Next, we estimate $\,I_3\,$, where we discuss each of the four terms occurring in the
definition of $Q^h$ (cf. \eqref{taylor2}) individually. In the following, we repeatedly omit the
time argument \betti{inside the integrals} for the sake of a shorter exposition.
At first, we have for every $\gamma>0$ (to be chosen later) that
\begin{align}\label{3.48}
&\left|\int_0^t\langle P(\bphi)(\chi-\zeta),A^{-1}(\psi-\chi)\rangle_V\ds\right|\\[1mm]
&\nonumber\le\,C_{16}
\int_0^t (\|\chi\|_{V'}+\|\zeta\|_{V'})\,(\|\psi\|_{V'}+\|\chi\|_{V'})\ds\\[1mm]
&\nonumber\le \,
\gamma\int_ 0^t\|\zeta\|_{V'}^2\ds\,+\,\frac{C_{17}}\gamma\int_0^t\left(\|\chi\|_{V'}^2
+\|\psi\|^2_{V'}\right)\ds\,,
\end{align} 
where in the first inequality we have used \eqref{ssbounds2} and 
the fact that (recall \eqref{gianni})
$$\|P(\bphi)(\chi-\zeta)\|_{V'}\,\le\,\|P(\bphi)\|_{W^{1,\infty}(\oma)}\,
\left(\|\chi\|_{V'}+\|\zeta\|_{V'} \right) \quad\mbox{a.\,e. in }\,(0,T). $$

\noindent
Next, \pier{in view of} \eqref{neu2}, we find that
\beqa\label{3.49}
&&\left|\int_0^t\langle (P(\vp^h)-P(\bphi))[(\s^h-\bs)-(\mu^h-\bmu)],A^{-1}(\psi-\chi)\rangle_V\,ds\right|\\[1mm]
&&\nonumber \le\,C_{18}\int_0^t\|(P(\vp^h)-P(\bphi))[(\s^h-\bs)-(\mu^h-\bmu)]\|_H\,(\|\psi\|_{V'}+\|\chi\|_{V'})\ds\\[1mm]
&&\nonumber \le \,C_{19}\,\|h\|^4_{\lzht}\,+\,C_{20}\int_0^t\left(\|\psi\|^2_{V'}\,+\,\|\chi\|_{V'}^2\right)\ds\,.
\eeqa

\noindent
In addition, \eqref{neu4} yields that
\begin{align}\label{3.50}
&\left|\int_0^t\langle P'(\bphi)\,(\bs-\bmu)\,\psi,A^{-1}(\psi-\chi)\rangle_V\ds\right|\\[1mm]
&\nonumber\le\,C_{21}\int_0^t\|P'(\bphi)(\bs-\bmu)\psi\|_H\,(\|\psi\|_{V'}+\|\chi\|_{V'})\ds\\[1mm]
&\nonumber\le\,\gamma\int_0^t\|\psi\|_V^2\ds\,+\,\frac{C_{22}}\gamma\int_0^t
\left(\|\psi\|^2_{V'}\,+\,\|\chi\|_{V'}^2\right)\ds\,.
\end{align}

\noindent
Finally, we obtain from \eqref{neu3} that
\begin{align}\label{3.51}
&\left| \int_0^t\langle (\bs-\bmu)\,R_2^h\,(\vp^h-\bphi)^2,A^{-1}(\psi-\chi)\rangle_V\ds\right|\\[1mm]
&\nonumber\le\,C_{23}\int_0^t\|(\bs-\bmu)\,R_2^h\,(\vp^h-\bphi)^2\|_{H}\,(\|\psi\|_{V'}+\|\chi\|_{V'})\ds \\[1mm]
&\nonumber\le\,
C_{24}\,\|h\|_{\lzht}^4\,+\,C_{25}\int_0^t\left(\|\psi\|^2_{V'}\,+\,\|\chi\|_{V'}^2\right)\ds\,.
\end{align}

Combining the estimates \eqref{3.45}--\eqref{3.51}, we have shown that for every $t\in [0,T]$ and $\gamma>0$
it holds that
\begin{align}\label{3.52}
&\frac 12\left(\|\psi(t)\|^2_{V'}\,+\,\|\chi(t)\|^2_{V'}\right)\,+
\,(1-\gamma)\int_0^t \|
\psi(s)\|_V^2\ds \,+\,\txinto |\chi|^2\dx\ds\\[1mm]
&\nonumber\le\,C_{26}\,\|h\|^4_{\lzht}\,+\,C_{27}\,\left(1+\gamma^{-1}\right)
\int_0^t\left(\|\psi(s)\|^2_{V'}\,+\,\|\chi(s)\|_{V'}^2\right)\ds\\[1mm]
&\nonumber\quad +\,C_{28}\,\|\psi\|^2_{\lzht}\,
 +\,\gamma\int_0^t\|\zeta(s)\|^2_{V'}\ds\,.
\end{align}
We still need to control the last two terms on the right-hand side of \eqref{3.52}. At first, notice that from \pier{the second equation in \eqref{3.44} and \eqref{neu1}  we can infer that}
\begin{align}\label{3.53}
&\int_0^t\|\zeta(s)\|_{V'}^2\ds\,\le\,C_{29}
\int_0^t\left(\|\psi(s)\|_V^2\,+\,\|\psi(s)\|_{V'}^2\right)\ds\\[1mm]
&\nonumber\quad \,+\,C_{30}
\int_0^t\|R_1^h(s)\,(\vp^h(s)-\bphi(s))^2\|_{V'}^2\ds\\[1mm]
&\nonumber\le\,C_{31}\int_0^t\left(\|\psi\|_V^2\,+\,\|\psi\|_{V'}^2\right)\ds\,+\,C_{32}\,\|h\|^4_{\lzht}\,.
\end{align}

\noindent Now observe that the compactness of the embeddings $V\subset H\subset V'$ implies that for every
$\gamma>0$ there is some constant $C_\gamma>0$ such that
\beq\label{ehren}
\|v\|_H^2\,\le\,\gamma\,\|v\|_V^2\,+\,C_\gamma\,\|v\|_{V'}^2\quad\forall\,v\in V.
\eeq
Hence, adjusting $\gamma>0$ appropriately small, 
invoking the estimates \eqref{3.52}--\eqref{ehren}, and applying Gronwall's lemma, we can finally infer that
\begin{align}\label{esti1}
&\|\psi^h\|_{L^\infty(0,t;V')\cap \lzvt}^2\,+\,\|\chi^h\|_{L^\infty(0,t;V')\cap\lzht}^2\,
+\,\|\zeta^h\|^2_{L^2(0,t;V')}\\[1mm]
&\nonumber\le C_{33}\,\|h\|^4_{\lzht}
\quad\forall\,t\in [0,T]\,.
\end{align}
 
\vspace{3mm}\noindent
\underline{\sc Second estimate:} \,\,\,At first, we observe that \eqref{gianni},
\eqref{esti1} and ({\bf H4}) imply that 
\beq\label{3.56}
\int_0^t\|P(\bphi(s))\,(\chi(s)-\zeta(s))\|_{V'}^2\ds\,\le\,C_{34}\,\|h\|^4_{\lzht}
\quad\forall\,t\in [0,T].
\eeq

\noindent Hence, it follows from \pier{\eqref{taylor2}, \eqref{neu2}--\eqref{neu4}, and} \eqref{esti1} 
that
\begin{align}\label{3.57}
&\|Q^h\|^2_{L^2(0,t;V')}\,\le\,C_{35}\,\|h\|^4_{\lzht}\quad\forall\,t\in [0,T].
\end{align} 

\noindent
Hence, testing \eqref{fs3} by $\chi^h$, we obtain, for every $t\in [0,T]$,
\begin{align}\label{3.58}
&\frac 12 \|\chi^h(t)\|_H^2\,+\,\txinto|\nabla\chi^h|^2\dx\ds
\,=\,-\txinto Q^h\,\chi^h\dx\ds\\[1mm]
&\nonumber\le \int_0^t\|Q^h(s)\|_{V'}\,\|\chi^h(s)\|_V\ds\,\le\,\frac 12\int_0^t\|\chi^h(s)\|_V^2\ds
\,+\,\frac 12\int_0^t\|Q^h(s)\|^2_{V'}\ds\,,
\end{align}
and Gronwall's lemma shows that
\beq\label{esti2}
\|\chi^h\|^2_{L^\infty(0,t;H)\cap \lzvt}\,\le\,C_{36}\,\|h\|_{\lzht}^4
\quad\forall\,t\in [0,T].
\eeq

\vspace{3mm}\noindent
\underline{\sc Third estimate:} \,\,\,We now multiply \eqref{fs2} by $\Delta\psi$ and
\betti{take $v=\psi$} in \eqref{fs1}, and integrate. Adding the resulting identities,
applying integration by parts and Young's inequality, and invoking the estimates \eqref{neu1}, \eqref{esti1}, and \eqref{3.57},
we easily obtain that, for every $t\in [0,T]$,
\begin{align}\label{3.60}
&\frac 12\,\|\psi(t)\|^2_H\,+\txinto|\Delta\psi|^2\dx\ds\\[1mm]
&\nonumber 
\,\le\txinto|\Delta\psi|\left(|F''(\bphi)\psi|+|R_1^h(\vp^h-\bphi)^2|\right)\dx\ds
\,+\int_0^t\|Q^h(s)\|_{V'}\,\|\psi(s)\|_V\ds\\[1mm]
&\nonumber
\,\le\,\frac 12\txinto|\Delta\psi|^2\dx\ds\,+\,C_{37}\,\|h\|^4_{\lzht}\,.
\end{align}
Using this and \eqref{fs2}, we have thus shown the estimate
\beq\label{esti3}
\|\psi^h\|_{L^\infty(0,t;H)\cap L^2(0,t;\hzwei)}^2 \,+\,\|\zeta^h\|_{\lzht}^2
\,\le\,C_{38}\,\|h\|^4_{\lzht}\quad\forall\,t\in [0,T].
\eeq

\noindent Comparison in \eqref{fs1} then yields that also
\beq\label{esti4}
\|\pt\psi^h\|_{L^2(0,t;W')}^2\,\le\,C_{39}\,\|h\|^4_{\lzht}
\quad\forall\,t\in [0,T].
\eeq

\vspace{3mm}\noindent
\underline{\sc Fourth estimate:} \,\,\,Now that \eqref{esti3} is shown, we also
have
\beq\label{3.63}
\int_0^t\|P(\bphi(s))(\chi(s)-\zeta(s))\|_H^2\ds\,\le\,C_{40}\,\|h\|^4_{\lzht}
\quad\forall\,t\in [0,T],
\eeq
which, together with \eqref{neu2}--\eqref{neu4}, implies the bound
\beq\label{3.64}
\|Q^h\|_{\lzht}^2\,\le\,C_{41}\,\|h\|_{\lzht}^4\quad\forall\,t\in [0,T].
\eeq

\noindent It is then an easy task (test \eqref{fs3} first by $\pt\chi^h$ and then
by $\,-\Delta\chi^h$) to see that also
\beq\label{esti5}
\|\chi^h\|_{H^1(0,t;H)\cap L^2(0,t;\hzwei)}^2\,\le\,C_{42}\,\|h\|^4_{\lzht}
\quad\forall\,t\in [0,T].
\eeq

\noindent With this, the inequality \eqref{3.30} is shown if we choose 
the function $Z(\lambda)$ as an appropriate multiple of $\lambda^4$. The
assertion is thus proved.\qed

\vspace{4mm}
\noindent
{\sc Remark \pier{4}} \quad Since the embedding of $H^1(0,T;W')\cap L^2(0,T;W)$
into $C^0([0,T];H)$ is continuous, we infer from Theorem 3.2 that the
control-to-state mapping ${\cal S}$ is also Fr\'echet differentiable into
$C^0([0,T];H)$ with respect to the first variable. From this it follows
that the reduced cost functional $\widetilde{\cal J}(u):={\cal J}
({\cal S}_1(u),u)$ (where ${\cal S}_1(u)$ denotes the first component of
${\cal S}(u)$) is Fr\'echet differentiable in $\ur$. Recalling that $\uad$ is a closed
and convex subset of $\lzq$, we conclude from standard arguments (which need no
repetition here) the
following result.

\vspace{4mm}\noindent
{\sc Corollary 3.3} \quad {\em Suppose that the assumptions $\as$ are fulfilled, and assume
that $\bu\in\uad$ is an optimal control for the problem $\PP$ with associated state
$(\bphi,\bmu,\bs)=\cs(\bu)$. Then we have}
\begin{align}
\label{vug1}
&\beta_Q\int_0^T\!\!\ioma (\bphi-\vp_Q)\,\xi \dx\dt\,+\,\beta_\oma\ioma(\bphi(T)-\vp_\oma)\,\xi(T)\dx\\[1mm]
&\nonumber \,\,+\,\beta_u\int_0^T\!\!\ioma \bu(v-\bu)\dx\dx\,\ge\,0 \quad\forall\,v\in\uad,
\end{align}
{\em where $\xi$ is the first component of the solution to the linearized system}
\eqref{ls1}--\eqref{ls5} {\em for $h=v-\bu$.}

\section{The control problem}
\label{control}
\setcounter{equation}{0}
\paragraph{Existence.} Consider the control problem $\PP$. We begin with the following existence result.

\vspace{3mm}\noindent
{\sc Theorem 4.1} \quad{\em Suppose that the assumptions $\as$ are fulfilled. Then the 
optimal control problem $\PP$ has a solution $\bu\in\uad$.}

\noindent
{\sc Proof:}  Let $\{u_n\}\subset\uad$ be a minimizing sequence for {\bf(CP)},
and let $(\vp_n,\mu_n,\s_n)={\cal S}(u_n)$,
$n\in \NN$. Then it follows from \eqref{ssbounds1} and (\ref{ssbounds2}) that there 
exist $(\vp,\mu,\s)$ and $u\in\uad$ such that, 
possibly for a subsequence which is again indexed by $n$, we have 
\begin{align}\no
&u_n\to u \quad\mbox{weakly star in }\,\liq,\\[1mm]
\nonumber
&\vp_n\to \vp \quad\mbox{weakly star in }\,H^1(0,T;V)\cap L^\infty(0,T;H^3(\oma)),\\[1mm]
&\nonumber\Delta\vp_n\to\Delta\vp \quad\mbox{weakly in }\,L^2(0,T;W),\\[1mm]
&\nonumber\mu_n\to\mu \quad\mbox{weakly star in }\,L^\infty(0,T;\betti{V})
\cap L^2(0,T;W),\\[1mm]
&\nonumber\s_n\to \s\quad\mbox{\pier{weakly star} in }\,H^1(0,T;H)\betti{{}\cap L^\infty(0,T;V){}}\cap L^2(0,T;W)\,.
\end{align}
In addition,
by virtue of standard compactness \pier{results (cf., e.g., \cite[Sec. 8, Cor. 4]{simon})}, we have the strong convergence 
\begin{align}\no
&\vp_n\to\vp \quad\hbox{strongly in } \pier{C^0([0,T]; H^2(\Omega))}\,, 
\end{align}
which implies, in particular,
that
\begin{align}\no
&\vp_n\to\vp \quad\hbox{strongly in } C^0(\overline{Q})\,,
\end{align}
whence also
\begin{align}\no
&F'(\vp_n)\to F'(\vp) \quad\mbox{and }\,\,\,P(\vp_n)\to P(\vp), \quad\mbox{both strongly in }\,
C^0(\overline{Q})\,.
\end{align}
In summary, we can pass to the limit as $n\to\infty$ in (\ref{ss1})--(\ref{ss5}), 
written for $(\vp_n,\mu_n,\s_n)$, 
finding that $(\vp,\mu,\s)={\cal S}(u)$; i.e., the pair $((\vp,\mu,\s),u)$ is 
admissible for {\bf (CP)}. 
It then follows from the weak sequential lower semicontinuity properties of ${\cal J}$ that 
$((\vp,\mu,\s),u)$ is an optimal pair for {\bf (CP)}. \qed

\paragraph{The adjoint system and first order necessary optimality conditions.} In order to 
establish the necessary first order optimality conditions for {\bf (CP)}, we need to eliminate $\xi$ from inequality \eqref{vug1}. 
To this end, we introduce the {\em adjoint system} which formally reads as follows:
\beqa\label{ad1}
&&-\partial_t  p+\Delta q-F''(\bphi)\, q+P'(\bphi)(\bs-\bmu)( r- p)
=\beta_Q\,(\bphi-\vp_Q) \quad\mbox{in }\,Q,\\[1mm]
\label{ad2}
&&q-\Delta p + P(\bphi)( p- r)=0 \quad\mbox{in }\,Q,\\[1mm]
\label{ad3}
&&-\partial_t r-\Delta r + P(\bphi)( r- p)=0 \quad\mbox{in }\,Q,\\[1mm]
\label{ad4}
&&\dn p \pier{{}=\dn q{}}  =\dn r=0 \quad\mbox{on }\,\Sigma,\\[1mm]
\label{ad5} 
&& r(T)=0, \quad  p(T)=\beta_\oma\,(\bphi(T)-\vp_\Omega) \quad\mbox{in }\,\Omega\,.
\eeqa

\noindent Since the final value $p(T)$ only belongs to $L^2(\Omega)$, we can at best expect the regularity 
\[ p\in  H^1(0,T; W')\cap C^0([0,T]; H)\cap L^2(0,T; W)\,,\]
which entails that \eqref{ad1} has to be understood in a weak variational sense. More precisely,
we call $(p,q,r)$ a solution to the adjoint system
\eqref{ad1}--\eqref{ad5} if and only if the functions $(p,q,r)$ satisfy the following conditions:
\beqa
\label{reguad}
&&p\in H^1(0,T;W')\cap C^0([0,T];H)\cap L^2(0,T;W),\quad  q\in\lzq, \\
&&\nonumber r\in H^1(0,T;H)\cap C^0([0,T];V)\cap L^2(0,T;W), 
\eeqa
the equations \pier{\eqref{ad1}--\eqref{ad5}} are satisfied almost everywhere in their respective
domains, \pier{but \eqref{ad1} and the related boundary condition in \eqref{ad4} 
hold} true in the sense that
\beqa\label{4.7}
&&\langle -\partial_t p(t),v\rangle_W\,+\ioma q(t)\Delta v\dx\,-\ioma F''(\bphi(t))\,q(t)\,v\dx\\[1mm]
&&\nonumber +\ioma P'(\bphi(t))(\bs(t)-\bmu(t))\,(r(t)-p(t))\,v\dx \,=\,\ioma\beta_Q\,(\bphi(t)-\vp_Q(t))v\dx
\\[2mm]
&&\nonumber\mbox{for all $v\in W$ and almost every $t\in (0,T)$}. 
\eeqa

We have the following existence and uniqueness result.

\vspace{4mm}\noindent
{\sc Theorem 4.2} \quad {\em \pier{Assume that the hypotheses $\as$  hold. Then} the adjoint system {\rm (\ref{ad1})--(\ref{ad5})} has a unique solution 
in the sense formulated above.}

\vspace{2mm}\noindent
{\sc Proof:} As in the proof of Theorem 3.1, we apply a Faedo-Galerkin approximation using the 
family $\{w_j\}_{j\in\NN}\subset W$ as a Galerkin basis in $W$ and ${\cal W}_n$ as
approximating finite-dimensional spaces.
Let $n\in\NN$ be fixed. We look for functions
of the form
\begin{equation*} 
p_n(x,t):=\sumk a_k^n(t)w_k(x), \quad r_n(x,t):=\sumk b_k^n(t)w_k(x), \quad
q_n(x,t):=\sumk c_k^n(t)w_k(x), 
\end{equation*}
which satisfy for almost every  $t\in (0,T)$ the following approximating problem:
\beqa
\label{ngal1}
&&(-\pt p_n(t),v)_H \,+\,(q_n(t),\Delta v)_H \,+\,(P'(\bphi(t))(\bs(t)-\bmu(t))(r_n(t)-p_n(t)),v)_H
\\[1mm]
&&\nonumber
-(F''(\bphi(t))\,q_n(t),v)_H\,=\,(\beta_Q(\bphi(t)-\vp_Q(t)),v)_H \quad\forall\,v\in {\cal W}_n,\\[2mm]
\label{ngal2}
&&(q_n(t),v)_H\,=\,(\Delta p_n(t) + P(\bphi(t))(r_n(t)-p_n(t)),v)_H 
\quad\forall\,v\in {\cal W}_n,\\[2mm]
\label{ngal3}
&&
(-\pt r_n(t),v)_H\,+\,\bigl(-\Delta r_n(t) +P(\bphi(t))(r_n(t)-p_n(t)),v\bigr)_H\,=\,0 
\quad\forall\,v\in {\cal W}_n,\\[2mm]
\label{ngal4}
&&r_n(T)=0, \quad p_n(T)=\mathbb{P}_n (\beta_\oma(\bphi(T)-\vp_\oma)),
\eeqa
where $\mathbb{P}_n$ denotes the orthogonal projector in $H$ onto ${\cal W}_n$.

Arguing as in the proof of Theorem 3.1, we can again infer that the backward-in-time
initial value problem \eqref{ngal1}--\eqref{ngal4} has a unique solution triple
$(p_n,q_n,r_n)\in (W^{1,\infty}(0,T;{\cal W}_n))^3$.
    
We now aim to derive a number of a priori estimates for the approximations. To this end, we denote by $C_i$, $i\in\NN$, positive constants that may depend on the data of the problem, but 
not on $n\in\NN$.   

\noindent \underline{\sc A priori estimates:}  \,\,\,Let $t\in [0,T]$ be arbitrary. We insert
$v=p_n(t)$ in \eqref{ngal1}, $v= \pier{{}- \Delta p_n(t)\in {\cal W}_n{}}$ in \eqref{ngal2}, and $v=r_n(t)$ in \eqref{ngal3},
add the resulting equations and integrate over $[t,T]$. In view of \eqref{ngal4},
we find the identity
\begin{align}\label{4.12}
&\frac 12\,(\|p_n(t)\|_H^2+\|r_n(t)\|_H^2)\,+\int_t^T\!\!\!\ioma\left(\pier{|\Delta p_n|^2}+
|\nabla r_n|^2\right)\dx\ds\\[1mm]
&\nonumber=\,\frac 12\,\|\mathbb{P}_n(\beta_\oma(\bphi(T)-\vp_\oma))\|_H^2\,+\,\int_t^T\!\!\!\ioma
F''(\bphi)\,p_n\,q_n\dx\ds \\[1mm]
&\nonumber\quad \,+\int_t^T\!\!\!\ioma\beta_Q(\bphi-\vp_Q)p_n\dx\ds
\,+\int_t^T\!\!\!\ioma P(\bphi)(p_n-r_n)(r_n\pier{{}+\Delta p_n{}})\dx\ds\\[1mm]
&\nonumber\quad\,-
\int_t^T\!\!\!\ioma P'(\bphi)(\bs-\bmu)(r_n-p_n)p_n\dx\,.
\end{align}
Using Young's inequality, it is easily seen that the first four summands on the right-hand
side are  bounded by an expression of the form
\beq\label{4.13}
C_1\,+\,\betti{\frac 12}\int_t^T\!\!\!\ioma\pier{|\Delta p_n|^2}\dx\ds\,+\,C_2\int_t^T\!\!\!\ioma\left(p_n^2
+r_n^2\right)\dx\ds\,,
\eeq
while for the last one (which we denote by $I$) it follows from H\"older's inequality\pier{, \eqref{ssbounds2} and the continuous embedding $V\subset L^3(\Omega)$} that, for any $\gamma>0$,
\begin{eqnarray}\label{4.14}
|I|&\!\!\le\!\!&C_3\int_t^T\left(\|\bs(s)\|_6+\|\bmu(s)\|_6\right)\|p_n(s)\|_2
\left(\|p_n(s)\|_3+\|r_n(s)\|_3\right)\ds\\[1mm]
\nonumber
&\!\!\le\!\!&\gamma\int_t^T\left(\|p_n(s)\|_V^2+\|r_n(s)\|_V^2\right)\ds\,+\,\frac{C_4}\gamma
\int_t^T\!\!\!\ioma\betti{p_n^2}\dx\ds\,.
\end{eqnarray} 
\pier{Hence, applying standard elliptic estimates to $p_n$ and adjusting $\gamma>0$
appropriately small, we deduce from Gronwall's lemma backward in time that
\beq\label{apri1}
\|p_n\|_{L^\infty(0,T;H)\cap L^2(0,T;\hzwei)}\,+\,
\|r_n\|_{L^\infty(0,T;H)\cap L^2(0,T;V)}\,\le\,  C_5\,.
\eeq
Next, taking $v= q_n(t) $ in \eqref{ngal2} and integrating in time, by \eqref{apri1} it is straightforward 
to deduce that 
\beq\label{pier2}
\|q_n\|_{ L^2(0,T;H)}\,\le\,  C_6\,.
\eeq
}%
Moreover, it is an easy task (by first inserting $v=-\pt r_n(t)$  and then
$v=-\Delta r_n(t)$ \pier{in \eqref{ngal3}}) to show that also
\beq\label{apri2}
\|r_n\|_{H^1(0,T;\betti{H})\cap L^\infty(0,T;V)\cap L^2(0,T;\hzwei)}\,\le\,\pier{C_7}\,.
\eeq

\vspace{3mm}
\noindent
\underline{{\sc Conclusion of the proof:}} 
It follows from the a priori estimates that there are functions $(p,q,r)$ such
that, possibly only for some subsequence which is again indexed by $n$,
\begin{align*}
p_n&\to p \quad\mbox{weakly in }\, L^2(0,T;W),\\[1mm]
q_n&\to q \quad\mbox{weakly in }\,\lzq,\\[1mm]
r_n&\to r \quad\mbox{weakly in }\,H^1(0,T;H)\cap L^2(0,T;\pier{W})\,,
\end{align*}
and, by continuous embedding, also
\begin{align*}
r_n&\to r \quad\mbox{weakly in }\,C^0([0,T];V)\,.
\end{align*}
It is now a standard matter (cf. the conclusion of the proof of Theorem 3.1) to show that
the triple $(p,q,r)$ is in fact a solution to the linear system \eqref{ad1}--\eqref{ad5}
having the asserted properties. Also the uniqueness can easily be proved;
we can allow ourselves to leave the argument to the interested reader.
\qed

\vspace{5mm}
We are now in the position to eliminate $\,\xi\,$ from (\ref{vug1}). We have the
following result. 

\vspace{4mm}\noindent
{\sc Theorem 4.3}  \quad
{\em Assume that the hypotheses $\betti{\as}$ are fulfilled, and suppose that
$\bu\in\uad$ is an optimal control for problem {\bf (CP)} with associated state 
$(\bphi,\bmu,\bs)={\cal S}(\bu)$ and adjoint state $(p,q,r)$. Then we have}
\begin{equation}\label{vug2}
\int_0^T\!\!\!\int_\Omega \left(r\,+\,\beta_u\,\bu\right)(v-\bu)\dx\dt\,\ge\,0
\quad\forall\,v\in\uad.
\end{equation}

\vspace{3mm}\noindent
{\sc Proof:}  We have, owing to \pier{\eqref{ls2}--\eqref{ls4}, \eqref{ss1weak}, 
\eqref{ad2}--\eqref{ad5}, and \eqref{4.7},} the following identities:
\begin{align}\label{id1}
0&=\int_0^T\!\!\!\ioma q\,[\eta+\Delta\xi-F''(\bphi)\xi]\dx\dt\,=\,
\int_0^T\!\!\!\ioma (q\,\eta-F''(\bphi)\,q\,\xi)\dx\dt\\[1mm]
\nonumber&\quad +\int_0^T\langle\pt p(t),\xi(t)\rangle_W\dt\,+\,\int_0^T\!\!\!\ioma
F''(\bphi)\,q\,\xi\dx\dt\\[1mm]
\nonumber&\quad -\int_0^T\!\!\!\ioma P'(\bphi)(\bs-\bmu)\,\xi\,(r-p)\dx\dt\,+
\int_0^T\!\!\!\ioma \beta_Q(\bphi-\vp_Q)\,\xi\dx\dt\,,
\end{align}

\vspace{-3mm}
\begin{align}
\label{id2}
0&=\int_0^T\langle\pt\xi(t),p(t)\rangle_V\dt\, +\int_0^T\!\!\!\ioma\nabla\eta\cdot
\nabla p \dx\dt 
\,-\int_0^T\!\!\!\ioma P(\bphi)(\rho-\eta)\,p\dx\dt\\[1mm]
\nonumber&\quad -\int_0^T\!\!\!\ioma P'(\bphi)(\bs-\bmu)\,\xi\,p\dx\dt\,,
\end{align}

\vspace{-3mm}
\begin{align}
\label{id3}
0&=\int_0^T\!\!\!\ioma r\,[\pt\rho-\Delta\rho +P'(\bphi)(\bs-\bmu)\,\xi
+P(\bphi)(\rho-\eta)-h]\dx\dt\\[1mm]
\nonumber &=\int_0^T\!\!\!\ioma \rho\,[-\pt r-\Delta r+P(\bphi)\,r]\dx\dt\,-
\int_0^T\!\!\!\ioma r\,h\dx\dt\\[1mm]
\nonumber&\quad +\int_0^T\!\!\!\ioma \xi\,P'(\bphi)(\bs-\bmu)\,r\dx\dt
\, -\int_0^T\!\!\!\ioma \eta\,P(\bphi)\,r\dx\dt\,. 
\end{align}
 Next, we employ integration by parts with respect to time \betti{in the second integral on the \pier{right-hand side of} \eqref{id1}} (which is permitted since
$p,\xi\in H^1(0,T;W')\cap L^2(0,T;W)$) to conclude that
\beq\label{id4}
I\,:=\int_0^T\langle \pt p(t),\xi(t)\rangle_W\dt\,=\,\ioma p(T)\,\xi(T)\dx
\,-\int_0^T\langle \pt\xi(t),p(t)\rangle_W\dt\,.
\eeq
Now observe that $p\in L^2(0,T;V)$ and $\pt\xi\in L^2(0,T;V')$. It then follows
from the second condition in \eqref{ad5} that
\beq\label{id5}
I\,= \ioma \beta_\oma(\bphi(T)-\vp_\oma)
\,\xi(T)\dx\,
-\int_0^T\langle \pt\xi(t),p(t)\rangle_V\dt\,.
\eeq 
Therefore, addition of the three identities \eqref{id1}--\eqref{id3} results in
\beq\label{id6}
0\,=\int_0^T\!\!\!\ioma \beta_Q(\bphi-\vp_Q)\,\xi\dx\dt\,+\,
\ioma \beta_\oma(\bphi(T)-\vp_\oma) \,\xi(T)\dx\,
-\int_0^T\!\!\!\ioma r\,h\dx\dt\,,
\eeq
and insertion of this identity in \eqref{vug1} yields the assertion.\qed

\vspace{4mm}
\noindent
{\sc Remark \pier{5}}  \quad The state system (\ref{ss1})--(\ref{ss5}), written for $(\vp,\mu,\s)
=(\bphi,\bmu,\bs)$, the adjoint system
and the variational inequality (\ref{vug2}) together form the first-order necessary optimality
conditions. Moreover, since $\uad$ is a nonempty, closed, and convex subset of $L^2(Q)$, 
\eqref{vug2} implies 
that for $\beta_u>0$ the optimal control $\bu$ is the $L^2(Q)$-orthogonal projection of 
$-\beta_u^{-1}\,r$ onto $\uad$, that is, we have
$$\bu(x,t)\,=\,\max\,\{u_{\rm min}(x,t),\,\min\,\{-\beta_u^{-1}r(x,t), 
u_{\rm max}(x,t)\}\,\} \quad\mbox{for a.\,e. }(x,t)\in Q. $$  


\end{document}